\documentclass[11pt]{article}
\usepackage[hscale=0.8,vscale=0.8]{geometry}

\usepackage{amsmath}
\usepackage{amsfonts}
\usepackage{amssymb}
\usepackage{amsthm}
\usepackage{mathtools}
\usepackage{algorithm}
\usepackage{algpseudocode}
\usepackage{graphicx}
\usepackage{caption}
\usepackage{subcaption}
\usepackage{multicol}
\usepackage{multirow}
\usepackage{hhline}
\usepackage{float}
\usepackage{color}
\usepackage{xcolor}
\usepackage{chngcntr}
\usepackage{tikz}
\usepackage{array}
\usepackage{comment}
\usepackage{geometry}
\usepackage{bm}
\usepackage{epsfig}

\usepackage{latexsym}
\usepackage{cite}
\usepackage{enumerate}
\usepackage{hyperref}
\usepackage{cleveref}

\graphicspath{{./figs/}}

\newcommand{\AG}[1]{{\color{red}#1}}
\newcommand{\yl}[1]{{\color{cyan}#1}}
\newcommand{\TL}{\textcolor{orange}}
\newcommand{\HELP}{\textcolor{blue}}
\usepackage[normalem]{ulem}

\def\wh{\widehat}
\def\wt{\widetilde}
\newcommand{\ds}{\displaystyle}

\def\qan{{\quad\hbox{and}\quad}}

\newcommand{\0}{{\mathbf{0}}}

\newcommand{\m}{{\mathbf{m}}}

\newcommand{\bv}{{\mathbf{v}}}
\newcommand{\by}{{\mathbf{y}}}

\def\bI{\mathbf{I}}

\newcommand{\bbP}{\mathbb{P}}
\newcommand{\bbR}{\mathbb{R}}

\newcommand{\bbZ}{\mathbb{Z}}

\newcommand{\bbeta}{{\boldsymbol\eta}}
\newcommand{\bzeta}{{\boldsymbol\zeta}}
\newcommand{\bPsi}{{\boldsymbol\Psi}}

\newcommand{\cN}{\mathcal{N}}
\newcommand{\cT}{\mathcal{T}}

\numberwithin{equation}{section}
\allowdisplaybreaks

\allowdisplaybreaks

\newtheorem{remark}{Remark}[section]

\title{A Structurally Informed Data Assimilation Approach for Nonlinear Partial Differential Equations}

\author{
  {\sc Tongtong Li}\thanks{Department of Mathematics, Dartmouth College, Hanover, NH 03755, USA, email: {\tt \{tongtong.li@dartmouth.edu, annegelb@math.dartmouth.edu, yoonsang.lee@dartmouth.edu\} }. All authors are supported in part by grant DoD ONR MURI \#N00014-20-1-2595. AG is also supported in part by grants NSF  DMS \#1912685, DOE ASCR  \#DE-ACO5-000R22725, and   AFOSR  \#FA9550-22-1-0411. YL is also supported by grant NSF DMS \#1912999.}	
\quad
    {\sc Anne Gelb}\footnotemark[1]~
\quad    
    {\sc Yoonsang Lee}\footnotemark[1]~}

\date{\today}

\begin{document}

\maketitle

\begin{abstract}

Ensemble-based Kalman filtering data assimilation is often used to combine available observations with numerical simulations to obtain statistically accurate and reliable state representations in dynamical systems.  However, it is well known that the commonly used Gaussian distribution assumption introduces biases for state variables that admit discontinuous profiles, which are prevalent in nonlinear partial differential equations. This investigation designs a new structurally informed prior that exploits statistical information from the simulated state variables. In particular, based on the second moment information of the state variable gradient, we construct a new weighting  matrix for the numerical simulation contribution in the data assimilation objective function.  This replaces the typical prior covariance matrix used for this purpose.  We further adapt our weighting matrix to include information in discontinuity regions via a clustering technique. Our numerical experiments demonstrate that this new approach yields more accurate estimates than those obtained using standard ensemble-based Kalman filtering on shallow water equations, even when it is enhanced with inflation and localization techniques.

\end{abstract}

\section{Introduction}\label{sec:introduction}
Data assimilation is a scientific process that combines observational data with mathematical models to obtain statistically accurate and reliable state representations for dynamical systems. It is widely used in a variety of scientific domains, notably atmospheric sciences, geoscience, and oceanographic sciences, and includes applications such as climate modeling, weather forecasting, and sea ice dynamics. Data assimilation techniques are often described to fit into one of two categories: (1) {\em filtering}, where the state estimation is continuously updated as new observations become available over time; and (2) {\em smoothing}, which re-evaluates past states using both historical and future observations so that previous estimations are retroactively improved. While both approaches have advantages and disadvantages, our goal here is to develop a new filtering method that synchronizes a given model with observation data. 

From a probabilistic perspective, (discrete-time) filtering is designed to sequentially update the probability distribution of the state variable conditioned on the accumulated data up to the current time. This update is realized in a two-step procedure. First, the forecast (prediction) step updates the probability distribution from the current to the next time step using a pushforward procedure driven by the given dynamic model, yielding a prior distribution. Bayes' formula is then applied in the analysis step to assimilate the observation data as the likelihood distribution, providing the posterior distribution. 

When the dynamical model is linear with additive, state-independent Gaussian model and observation errors, it is known that the Kalman filter \cite{Kalman60} provides the optimal data assimilation method \cite{Cohn97}, in the sense that it determines a complete characterization of the posterior distribution in terms of its mean and covariance. Problems of interest generally depart from these assumptions, however. In this regard there are a variety of partial solutions that fit the more realistic nonlinear and/or non-Gaussian cases, each better suited for particular environments. For example, sequential Monte Carlo algorithms such as particle filters \cite{Gordon93} provide robust solutions for low dimensional states and computationally inexpensive likelihoods. These techniques use either a reweighting or resampling step to convert the prior ensemble to a posterior ensemble after each observation. Unfortunately the weights degenerate for high-dimensional problems. In particular, the huge state and observation dimensions, along with the highly nonlinear evolution operators, make these methods impractical for geophysical problems \cite{Snyder08}. Alternatives include the extended Kalman filters based on linearizations of nonlinear dynamics and other Gaussian filters like the unscented Kalman filter \cite{Julier00}.  Such approaches involve extensive use of numerical approximations however, and in general provide no guarantees for accurate or meaningful uncertainty quantification \cite{Law12}. 

In this more realistic environment, robust results can be obtained using ensemble-based Kalman filters \cite{Evensen94}. Rather than explicitly engaging the entire state distribution as the standard Kalman filter does, ensemble-based Kalman filters store, propagate, and update an ensemble of vectors that collectively approximate the state distribution. Specifically, in the forecast step a particle approximation of the filtering distribution is propagated through the transition dynamics to yield a forecast ensemble at the next observation time. This ensemble is then updated in the analysis step to estimate the posterior distribution. For example, ensemble square-root Kalman filters (ESRF) \cite{Tippett03} convert the prior ensemble to a posterior ensemble through a linear transformation. Although effective in reducing computational cost, the covariance matrix, which plays an essential role from the prior to posterior shift,  is often poorly approximated.

All that said, the ensemble-based Kalman filters do not always provide the most suitable characterization of the posterior distribution for non-Gaussian models \cite{Mandel11}. These difficulties become more pervasive in geophysical applications  such as sea ice modeling, where the sea ice state variable often contains sharp features in the sea ice cover, such as leads and ridges. It is therefore imperative that assimilation schemes used for such problems are able to preserve the physical principles of the underlying problem while also accurately recovering solutions with discontinuous profiles, even when such features are not directly observed. Methods have been developed in this direction. For example, in \cite{Srang14}, when estimating discontinuous systems,  the technique deliberately chooses basis functions to  include discontinuities for the unscented transform in the continuous-discrete unscented Kalman filter (CDUKF). In \cite{Levy10} a decision criteria from a fuzzy verification metric is adopted to update the model state. The method assumes that there exist both adequate model validation and verification measures for lower dimensional features such as discontinuities. 

Our investigation is partially inspired by the equivalency of these aforementioned statistical data assimilation techniques \cite{Kalman60, Julier00, Evensen94, Tippett03, Srang14} to various forms of variational methods \cite{Kalnay02}. Specifically, since the filtering distribution in Kalman filtering is Gaussian, then by Bayes' formula the posterior mean is equivalent to the maximum a posterior (MAP) estimate, and can therefore be obtained by minimizing the negative logarithm of the posterior distribution. From the variational approach perspective, this is equivalent to solving an inverse problem where the objective function is the misfit between the observations and model state coupled with a $\ell_2$ regularization term defined by the difference between the optimal state and prior solution by forecast model. The relationship between the variational approach and Tikhonov regularization was described in \cite{Johnson05} by reformulating the objective function and using singular value decomposition techniques. It is also possible to replace the $\ell_2$-norm regularization  by another norm when appropriate. For example, $\ell_1$ regularization methods, often referred to as compressive sensing for sparse signal recovery \cite{CT}, are extensively used for image recovery since an image is typically sparse in the edge domain. The $\ell_1$ penalty approach has also been adopted under the variational data assimilation framework \cite{Freitag10}, and when the state variable admits sparse representation in a specific (e.g.~gradient) domain, a mixed variational regularization approach was shown in some instances to yield better results \cite{Freitag13, Asadi19}. 

While heuristically the $\ell_1$ penalty can be viewed as a term encoding a prior belief quantifying an assumption regarding the state variable (e.g.~sparse gradient domain), determining an appropriate prior in general is notoriously difficult. Indeed, the standard $\ell_1$ prior does not provide any {\em local} information regarding the sparse domain.  This was addressed in the context of image recovery in \cite{Adcock19,GelbScarnati2019}, where a weighted $\ell_1$ joint sparsity algorithm was developed based on variance from multiple measurement information.\footnote{Such observation data is often referred to as {\em multiple measurement vectors} (MMVs), where each measurement vector typically contains incomplete information regarding the underlying image of interest.  Their mutual (joint) information is exploited to improve the recovery.}  This variance-based  weighted $\ell_1$ penalty term  was then incorporated into the Bayesian framework as a {\em support informed} prior in \cite{Gelb22}.  Since an ensemble can be viewed similarly to  MMV data, here we are inspired to develop a similar idea as was used in \cite{Gelb22} that will more accurately combine the observation data into the data assimilation update, as we will now describe. 

The ensembles in the ensemble-based Kalman filters play a natural role in providing such empirical statistics. In particular, updates from the prior to the posterior may be understood in the context of regularization, with the model/data fidelity balanced by the covariance of the prior ensembles. Since the prior sample covariance does not capture local discontinuity structure, it is not an optimal weighting choice.  

In this paper we consider the ensemble-based Kalman filter known as the {\em ensemble transform Kalman filter} (ETKF) \cite{Bishop01}, which is a variant of the ESRF. Inspired by the support informed prior idea in \cite{Gelb22} we reformulate the analysis step using a regularization framework that provides local information. Specifically, we design our prior by constructing a new weighting matrix based on the second moment information of the state variable gradient, noting that more variability occurs in the neighborhoods of discontinuities, that is, within the support of the state variable in the gradient space. Moreover, since the prior information is weighted more heavily near the discontinuity region, we see that  accurate state variable approximation is also closely linked to the choice of prediction method in the forecast step. Because of this, it is important to carefully choose the numerical scheme used to solve the dynamic model, and more specifically, the method must be able to capture discontinuity information. The  weighted essentially non-oscillatory scheme (WENO) \cite{Liu94}, which is known to yield higher-order convergence for smooth solutions while resolve the discontinuities in the sharp feature in a non-oscillatory manner, provides one such method. Finally, we further adapt the weighting matrix to include neighborhood information surrounding discontinuities when updating the state solution in regions containing discontinuities. We refer to this as {\em clustering}.

We test our new gradient second moment weighting strategy on the dam break problem governed by the shallow water equations. The results demonstrate a more accurate approximation than currently used methodology. We also note that since our new weighting matrix is constructed from the ensemble, which is a low rank approximation, our approach does not add computational complexity. Finally, as the minimization is still based on the $\ell_2$ norm, there is a closed-form solution so that no extra optimization algorithms are needed, as would be the case for the $\ell_1$ or mixed $\ell_1$ and $\ell_2$ prior. 

The rest of the paper is organized as follows. In Section \ref{sec:preliminaries} we review the standard data assimilation techniques and its relation to minimization approach. We propose our method in Section \ref{sec:method}, and provide some numerical examples in Section \ref{sec:numerical}.  Section \ref{sec:conclusion} provides some concluding remarks.

\section{Preliminaries}
\label{sec:preliminaries}
We begin by establishing the notation  used throughout our manuscript, which mainly follows \cite{Law15}. Let $\bbZ^{+}$ denote the set of non-negative integers and $C(X,Y)$ denote the set of continuous functions mapping from space $X$ to space $Y$. We also use $\0$ to represent the null vector or tensor as appropriate.  The Hadamard (entrywise) product of matrices $M:=(M_{i,i'})_{i,i'=1,\cdots,n}$ and $N=(N_{i,i'})_{i,'=1,\cdots,n}$ is given by
$$M \odot N := (M_{i,i'}N_{i,i'})_{i,i'=1,\cdots,n}.$$
Finally, given vector $\bv$ and a symmetric positive definite matrix $C$, we denote the corresponding weighted norm (squared) as
\begin{equation}
\vert \bv \vert_C^2:= \vert C^{-\frac{1}{2}} \bv \vert^2 = \bv^T C^{-1} \bv.\label{eq:normsquared}
\end{equation}

\subsection{Problem set up}
\label{subsec:setup}
Let $v(x,t)$, for location $x\in \Omega$ and time $t \in [0,T]$,  be the variable of interest which determines the state of the system on domain $\Omega$ for given time $t$. We discretize the spatial domain and define its numerical approximation $\bv_j:=\bv(t_j) \in \bbR^n$ to be the corresponding vector at time instance $t_j = j\Delta t$, where $\Delta t = \frac{T}{J}$ and $J\in \bbZ^{+}$ is the final time step for simulation. We assume that the solution trajectory of $\bv$ follows the deterministic discrete model 
\begin{equation}\label{model: dynamic}
\bv_{j+1}=\bPsi(\bv_j), \quad j = 0,\dots, J-1,
\end{equation}
where the discrete operator $\bPsi: \bbR^n \rightarrow \bbR^n$ is often defined to approximate a known system of partial differential equations (PDEs), which is itself a low rank  approximation of the underlying dynamics. We further assume that the initial state $\bv_0$ follows a Gaussian distribution with mean $\m_0$ and covariance matrix $C_0$:
\begin{equation}\label{assump: initial Gaussian}
\bv_0 \sim \cN(\m_0, C_0). 
\end{equation}
Finally, we assume that the given observational data $\by \in \mathbb{R}^m$ are defined as 
\begin{equation}\label{model: observation}
\by_{j+1}=h(\bv_{j+1})+\bbeta_{j+1}, \quad j = 0,\dots, J-1,
\end{equation}
where $h: \bbR^n \rightarrow \bbR^m$ is the (possibly nonlinear) observation operator and each $\bbeta_{j+1}$ is an independent and identically distributed (i.i.d.) sequence, independent of $\bv_0$, with 
\begin{equation}\label{assump: noise Gaussian}
\bbeta_1 \sim \cN(\0,\Gamma).
\end{equation}
Here $\Gamma$ is assumed to be known and symmetric positive definite (SPD).

\subsection{Filtering}
\label{subsec:filter}
Filtering can be described as sequentially updating the probability distribution of the state as new data are acquired.  In essence the conditional probability density function $\bbP(\bv_{j+1} \mid \by_{1:j+1})$ is obtained by sequentially updating $\bbP(\bv_{j} \mid \by_{1:j})$ using $\by_{1:j}:=\{ \by_1, \cdots, \by_j\}$ through a two-step process. First in the {\em prediction step}, $\bbP(\bv_{j+1} \mid \by_{1:j})$ is calculated from $\bbP(\bv_{j} \mid \by_{1:j})$ using the Markov kernel
\begin{equation}\label{eq: pdf 1}
\bbP(\bv_{j+1} \mid \by_{1:j}) = \int \bbP(\bv_{j+1} \mid \by_{1:j}, \bv_{j}) \bbP(\bv_{j} \mid \by_{1:j}) d\bv_{j}, 
\end{equation}
through application of the dynamics operator $\bPsi$. This is followed by the {\em analysis step} in which  the {\em posterior} $\bbP(\bv_{j+1} \mid \by_{1:j+1})$ is obtained from the {\em prior} $\bbP(\bv_{j+1} \mid \by_{1:j})$ via Bayes' formula
\begin{equation}\label{eq: pdf 2}
\bbP(\bv_{j+1} \mid \by_{1:j+1}) \propto \bbP(\by_{j+1} \mid \bv_{j+1}) \, \bbP (\bv_{j+1} \mid \by_{1:j}), 
\end{equation}
where the \emph{likelihood} $\bbP(\by_{j+1} \mid \bv_{j+1})$ is determined from the observational data. Here $\propto$ implies that the two sides are equal to each other up to a multiplicative constant that does not depend on $\bv$ or $\by$.

\subsubsection{Kalman filtering (KF)}
\label{subsec:kalmanfilter}
If $\bPsi$ in \eqref{model: dynamic} and $h$ in \eqref{model: observation} are both linear, that is,
\begin{equation}
\label{eq:linearcase}
\bPsi(\bv)= M \bv, \ h(\bv)=H \bv,
\end{equation}
then  due to the Gaussian assumptions on initial state $\bv_0$ (c.f.~\eqref{assump: initial Gaussian})  and observation noise $\{\bbeta_{j+1} \}_{j\in \bbZ^{+}}$ (c.f.~\eqref{assump: noise Gaussian}), the filtering distribution posterior in \eqref{eq: pdf 2} must also be Gaussian.  Furthermore, it is entirely characterized by its mean and covariance as 
\begin{equation}
\label{eq:postdistributionlinear}
\bv_{j} \mid \by_{1:j} \sim  \cN(\m_j,C_j),
\end{equation}
where the posterior mean and covariance at time $t_j$  are respectively denoted by $\m_j$  and $C_j$. Similarly, the prior distribution at time $t_{j+1}$ is given by 
\begin{equation}
\label{eq:priordistributionlinear}
\bv_{j+1} \mid \by_{1:j} \sim  \cN(\wh{\m}_{j+1},\wh{C}_{j+1}),
\end{equation}
with prior mean and covariance respectively $\wh{\m}_{j+1}$ and $\wh{C}_{j+1}$.\footnote{Throughout our exposition, unless otherwise stated, we will use hat notation, e.g.~$(\wh{\m}, \wh{C})$, to denote estimates related to the prior, while estimates related to the posterior will not, e.g.~$(\m, C)$.} Based on \eqref{eq:linearcase}, the KF prediction step 
updates the prior mean $\wh{\m}_{j+1}$ as
\begin{equation}
\wh{\m}_{j+1} = M \, \m_j, \label{eq:hatmean_update}
\end{equation}
with corresponding covariance
\begin{equation}
\wh{C}^{KF}_{j+1}:= \wh{C}_{j+1} = M\, C_j\, M^T. \label{eq:KFcov_update}
\end{equation}
In the analysis step the posterior mean and covariance satisfy
\begin{eqnarray}
C_{j+1}^{-1}&=&\wh{C}_{j+1}^{-1} + H^T \Gamma^{-1} H,\nonumber\\
C_{j+1}^{-1} \m_{j+1} &=& \wh{C}_{j+1}^{-1} \wh{\m}_{j+1}+H^T \Gamma^{-1}\by_{j+1}. \label{eq: KF post 1}
\end{eqnarray}
Observe that \eqref{eq: KF post 1} is a direct consequence of Bayes rule and the Gaussianity of \eqref{eq: pdf 2}, that is
\begin{equation}\label{eq: KF post}
\exp \left( -\frac{1}{2} \vert \bv_{j+1} - \m_{j+1} \vert^2_{C_{j+1}} \right) \propto \exp \left(  -\frac{1}{2}\vert \by_{j+1}-H \bv_{j+1} \vert_{\Gamma}^2 -\frac{1}{2} \vert \bv_{j+1}-\wh{\m}_{j+1} \vert_{\wh{C}_{j+1}}^2 \right).
\end{equation}
Equivalently, the Kalman filter algorithm  updates the posterior mean and covariance through the iterative procedure 
\begin{subequations}
\label{eq: KF post 4}
\begin{equation}
\m_{j+1}=(I-K_{j+1}H)\wh{\m}_{j+1}+ K_{j+1}\by_{j+1}, \label{eq:mean_update}
\end{equation}
\begin{equation}
C_{j+1}=(I-K_{j+1}H) \wh{C}_{j+1}, \label{eq:cov_update}
\end{equation}
\end{subequations}
where
\begin{subequations}
\begin{equation}
K_{j+1}=\wh{C}_{j+1}H^{T}S_{j+1}^{-1}, \label{eq:gain_update}
\end{equation}
\begin{equation}
S_{j+1} =H \wh{C}_{j+1}H^{T} + \Gamma. \label{eq:S_update}
\end{equation}
\end{subequations}
Here $K_{j+1}$ is known as the {\em Kalman gain}. The equivalence of \eqref{eq: KF post 4} and  \eqref{eq: KF post 1} can be demonstrated using the Sherman-Morrison-Woodbury formula \cite{Sherman50, Woodbury50}. Since the matrix inversion in \eqref{eq: KF post 4} is typically performed in a lower dimensional data space, it is generally more computationally efficient.

\begin{remark}\label{rem:cov_update}
We note that \eqref{eq:cov_update} is also used when $\wh{C}_{j+1}$ is determined from constructions that may differ from the one given by \eqref{eq:KFcov_update}, as will be seen in Section \ref{sec:ETKF}.
\end{remark}
\subsubsection{Ensemble Kalman filtering (EnKF)}
\label{sec:ensemblefiltering}

Kalman filtering \cite{Kalman60} provides a complete solution for the sequential update of the probability distribution of the state given observational data \cite{Cohn97, Law15}.  The problems we are interested in often  do not satisfy the assumptions regarding linearity and Gaussian distributions, however.  Ensemble-based Kalman filters, for which an ensemble of particles is employed to inform the filtering update, are used to improve the performance outcomes in such non-idealistic cases \cite{Evensen94,Evensen06}. For context, the general ensemble-based Kalman filter methodology is briefly reviewed below.

The ensemble $\{ \bv^{(k)}_{j} \}_{k=1}^{K}$ is constructed by sampling from the filtering distribution $\bbP(\bv_{j} \mid \by_{1:j})$ 
at time $t_j$. The prediction step is generated by applying the dynamic model \eqref{model: dynamic} to each ensemble member yielding
\begin{equation}\label{eq: EnKF predict 1}
\wh{\bv}^{(k)}_{j+1}=\bPsi(\bv^{(k)}_j).
\end{equation}
We then compute the sample prior mean and covariance matrix as
\begin{subequations}
\begin{equation}
\ds \wh{\m}_{j+1}=\frac{1}{K}\sum_{k=1}^{K}\,\wh{\bv}^{(k)}_{j+1}, \label{eq: EnKF prior 1}
\end{equation}
\begin{equation}
\ds \wh{C}^{EnKF}_{j+1}:= \wh{C}_{j+1} =\frac{1}{K-1} \sum_{k=1}^{K}\,(\wh{\bv}^{(k)}_{j+1}-\wh{\m}_{j+1})(\wh{\bv}^{(k)}_{j+1}-\wh{\m}_{j+1})^{T},  \label{eq: EnKF prior 2}   
\end{equation}
\end{subequations}
noting that they respectively serve as low-rank representations of the prior mean and covariance.\footnote{For ease of presentation, as long as the context is clear, we use the same notation to denote both sample prior mean and prior mean, as well as for the sample covariance matrix and covariance matrix. Moreover, we replace the cumbersome notation  $\wh{C}^{EnKF}_{j+1}$ in \eqref{eq: EnKF prior 2} with $\wh{C}_{j+1}$ moving forward since our proposed method employs this definition.}
Ensemble Kalman filter (EnKF) \cite{Evensen94}, updates the prior ensembles $\{ \wh{\bv}^{(k)}_{j+1} \}_{k=1}^{K}$ from the new observation data $\by_{j+1}$ using a stochastic update step based on the perturbed observations. This is accomplished by generating a set of simulated observations $\{ \by^{(k)}_{j+1} \}_{k=1}^{K} = \by_{j+1} + \bbeta_{j+1}^{(k)}$, where  $\bbeta_{j+1}^{(k)} \sim \cN(0,\Gamma)$.  Each ensemble of particles is then updated according to 
\begin{eqnarray}
\bv_{j+1}^{(k)}&=&(I-K_{j+1}H)\wh{\bv}_{j+1}^{(k)}+ K_{j+1}\by_{j+1}^{(k)},
\label{eq:perturbedobservation}
\end{eqnarray}
with $K_{j+1}$ as defined in \eqref{eq:gain_update} using \eqref{eq: EnKF prior 2}.   Note that \eqref{eq:perturbedobservation}
can be equivalently understood as a family of minimization problems given by
\begin{gather*}
\bv_{j+1}^{(k)} = \underset{\bv}{\text{arg min}} 
\left(\frac{1}{2}\vert \by_{j+1}^{(k)}-H \bv \vert_{\Gamma}^2 + \frac{1}{2} \vert \bv-\wh{\bv}_{j+1}^{(k)} \vert_{\wh{C}_{j+1}}^2\right).
\end{gather*}
A disadvantage of EnKF is that, by design, \eqref{eq:perturbedobservation} introduces additional sampling errors. 

\subsubsection{Ensemble transform Kalman filtering (ETKF)}
\label{sec:ETKF}
The class of ensemble square-root Kalman filters (ESRFs) \cite{Tippett03} provides a deterministic update from prior to posterior ensembles, and therefore does not rely on such perturbed observations. Specifically, each prior ensemble member is shifted and scaled through a linear transformation process so that the posterior ensembles have mean and covariance satisfying the Kalman updates given by  \eqref{eq: KF post 4}. This is realized as 
\begin{equation}
\ds \bv^{(k)}_{j+1} = \m_{j+1}+\bzeta^{(k)}_{j+1},\quad k = 1,\dots,K,
\label{eq:zeta}
\end{equation}
where each $\bzeta^{(k)}_{j+1}$ is designed to follow the Gaussian distribution $\cN(\0,C_{j+1})$.  
As a result, only one minimization problem is needed to solve to obtain the posterior mean:
\begin{equation}
\m_{j+1} = \underset{\bv}{\text{arg min}} 
\left(\frac{1}{2}\vert \by_{j+1}-H \bv \vert_{\Gamma}^2 + \frac{1}{2} \vert \bv-\wh{\m}_{j+1}^{(k)} \vert_{\wh{C}_{j+1}}^2\right),
\label{eq:ESRF}
\end{equation}
followed by a transformation process \eqref{eq:zeta} to obtain the posterior ensembles.
In this regard there are many variants of ESRF, with the ensemble adjustment Kalman filter (EAKF) \cite{Anderson01} and the ensemble transform Kalman filter (ETKF) \cite{Bishop01} being the most widely used.  Since our new approach modifies the ETKF, we describe below how it is used to update the posterior ensembles for \eqref{eq:zeta}. We emphasize that our main contribution, as described in Section \ref{sec:method}, is in developing a new structurally informed prior to be generally used in \eqref{eq:ESRF}, regardless of implementation.

The main feature of the ETKF is that the updates of the ensemble particles are realized by finding an appropriate transformation operator. 
To this end we define the prior centered ensemble as
$$ \wh{X}_{j+1} = \frac{1}{\sqrt{K-1}} \left[ \wh{\bv}^{(1)}_{j+1}-\wh{\m}_{j+1}, \cdots, \wh{\bv}^{(K)}_{j+1}-\wh{\m}_{j+1} \right].$$
Substituting this expression into \eqref{eq: EnKF prior 2} yields the equivalent expression for the sample prior covariance
$$\wh{C}_{j+1}=\wh{X}_{j+1}\,\wh{X}_{j+1}^T.$$ 
We next define the posterior centered ensemble as
\begin{equation}\label{eq: ETKF X}
X_{j+1}=\wh{X}_{j+1} T^{\frac{1}{2}}_{j+1},
\end{equation}
where the transformation operator $T_{j+1}$ is given by
\begin{equation}\label{eq: ETKF T}
T_{j+1}=\left[ I+(H\wh{X}_{j+1})^T \,\Gamma^{-1}(H\wh{X}_{j+1}) \right]^{-1}.
\end{equation}
Incorporating \eqref{eq: ETKF T} into \eqref{eq:cov_update} yields the posterior covariance update
\begin{equation}\label{eq: ETKF C}
C_{j+1}=X_{j+1}\,X_{j+1}^T.
\end{equation}
Finally, since the perturbations $\{\bzeta^{(k)}_{j+1}\}_{k = 1}^K$ in \eqref{eq:zeta} are defined such that each follows the distribution $\cN(\0,C_{j+1})$, \eqref{eq: ETKF C} provides their explicit determination as
\begin{equation}\label{eq: ETKF zeta}
X_{j+1} = \frac{1}{\sqrt{K-1}} \left[ \bzeta^{(1)}_{j+1}, \cdots, \bzeta^{(K)}_{j+1} \right].
\end{equation}
\subsection{Inflation and localization}
\label{sec:inflatelocal}
For the practical implementation of ETKF, we briefly describe two important techniques commonly used in data assimilation. The first,  {\em variance inflation} \cite{Anderson07, Law15}, mitigates the phenomenon of filtering that causes an increasingly tight prior distribution which is decreasingly influenced by new observations for long time simulations.  This undesirable result is due to the biased estimates of the covariance matrix that may come from, for instance, a small ensemble size \cite{Furrer07}. For a given level of observational noise, the model variance scale can be ``inflated" to ensure that new data are adequately observed. Different inflation approaches, for example using additive or multiplicative constructions, will ultimately impact the value of $\wh{C}_{j+1}$. Here we adopt a simple form of multiplicative inflation, which in the context of ETKF is achieved by rescaling the centered ensemble as
\begin{equation}
\wh{X}_{j+1} = \alpha \wh{X}_{j+1}, \quad \alpha > 1,\label{eq:inflation}
\end{equation}
so that the prior ensemble is consistent with the inflated covariance.

The second technique commonly used is  {\em localization} \cite{Evensen06}. The idea here is to reduce unwanted correlations in $\wh{C}_{j+1}$ between state components that are separated by large distances in physical space. It is reasonable to assume that the correlation between state components generally decays as they get further away from each other, which also implies that this correlation is increasingly influenced by the sample error as the states are further apart. Spurious correlations can also cause rank deficiency, which means that traditional optimization methods will fail when solving the corresponding minimization problem(s). {\em Covariance localization} can be used to pre-process $\wh{C}_{j+1}$ and improve its conditioning by replacing 
\begin{equation}
\wh{C}_{j+1} = \wh{C}_{j+1} \odot \cT,\label{eq:localization}
\end{equation}
where $\cT$ is some sparse SPD matrix. Localization can be understood as a pre-processor that convolves the empirical correlation matrix with a given localized kernel, such as the compactly supported fifth-order piecewise polynomial correlation function \cite{Gaspari99}. In our numerical experiments we use a simplified version of localization, which we implement by extracting some prescribed bandwidth of the covariance matrix $\wh{C}_{j+1}$. Specifically, we define $\cT$ to be a binary Toeplitz matrix. In so doing, we avoid the additional complexity that comes with more parameter tuning while also improving computational efficiency.  We note that inflation and localization are typically used in conjunction with one another.  Algorithm \ref{alg:ETKF} describes the ETKF framework using these techniques.

\begin{algorithm}[!ht]
\caption{ETKF method with inflation and localization techniques}\label{alg: hybrid}
\begin{algorithmic}
\State Initialize ensembles $\bv^{(k)}_0$, $k = 1,\dots,K$. Choose inflation parameter $\alpha$ in \eqref{eq:inflation} and localization matrix $\cT$ in \eqref{eq:localization}.
\For{$j=1$ to $J-1$} 
\State Compute prior ensembles $\ds \wh{\bv}^{(k)}_{j+1}=\bPsi(\bv^{(k)}_j)$ for $k=1, \cdots, K.$ 
\If{data observations are available at time instance $j+1$}
\State Compute prior mean $\ds \wh{\m}_{j+1}=\frac{1}{K}\sum_{k=1}^{K}\,\wh{\bv}^{(k)}_{j+1}$.
\State Compute inflated centered ensemble $\ds \wh{X}_{j+1} = \frac{\alpha}{\sqrt{K-1}} \left[ \wh{\bv}^{(1)}_{j+1}-\wh{\m}_{j+1}, \cdots, \wh{\bv}^{(K)}_{j+1}-\wh{\m}_{j+1} \right]$ and covariance matrix $\ds \wh{C}_{j+1}=\wh{X}_{j+1}\,\wh{X}_{j+1}^T$.
\State Compute the transformation matrix $T_{j+1}$ and $X_{j+1}$ according to \eqref{eq: ETKF T} and \eqref{eq: ETKF X} respectively.
\State Solve the minimization problem for the posterior mean
\begin{gather*}
\m_{j+1} = \underset{\bv}{\text{arg min}} 
\left(\frac{1}{2}\vert \by_{j+1}-H \bv \vert_{\Gamma}^2 + \frac{1}{2} \vert \bv-\wh{\m}_{j+1} \vert_{\wh{C}_{j+1}\odot \cT}^2\right), 
\end{gather*}
\State Update the posterior ensembles according to the transformation \eqref{eq:zeta}
where $\bzeta^{(k)}_{j+1}$ is defined by \eqref{eq: ETKF zeta}.
\Else{ $\ds \bv^{(k)}_{j+1}=\wh{\bv}^{(k)}_{j+1}$ for $k=1, \cdots, K.$}
\EndIf
\EndFor
\end{algorithmic}
\label{alg:ETKF}
\end{algorithm}

\section{A new structurally informed prior design}
\label{sec:method}
What is clearly evident from Section \ref{sec:preliminaries} is the critical role played by prior covariance $\wh{C}_{j+1}$ in balancing the model and data information for the minimization problem \eqref{eq:ESRF}. Notably, if the underlying solution has a discontinuous profile, the prior covariance obtained using \eqref{eq: EnKF prior 2} no longer yields accurate results even after inflation and localization techniques are applied.  In particular, the construction of the covariance does not capture any information regarding discontinuities --  important features that by all accounts should be incorporated into the prior belief. Hence we propose to reweight the objective function in order to capture such information.
That is, we essentially seek a solution that appropriately resembles the information coming from both the observation and the prior. For simplicity and ease of notation, in our development we focus on a single  time instance yielding the minimization problem 
\begin{equation}
\label{eq:obj_w}
\m = \underset{\bv}{\text{argmin}}\left(\frac{1}{2}\vert \by-H \bv \vert_{\Gamma}^2 + \frac{1}{2} \vert \bv-\wh{\m} \vert_{W}^2\right), 
\end{equation}
where the weighted norms are defined by \eqref{eq:normsquared}.  Our goal is therefore to determine the weighting matrix $W$ in \eqref{eq:obj_w} that yields the best results. We consider two alternative scenarios: an ideal environment, where the sampled spatial observations are dense (Section \ref{sec:full}) and a less ideal environment, where they are sparse (Section  \ref{sec:sparse}).

\begin{figure}[h!]
\centering
\includegraphics[width=0.32\textwidth]{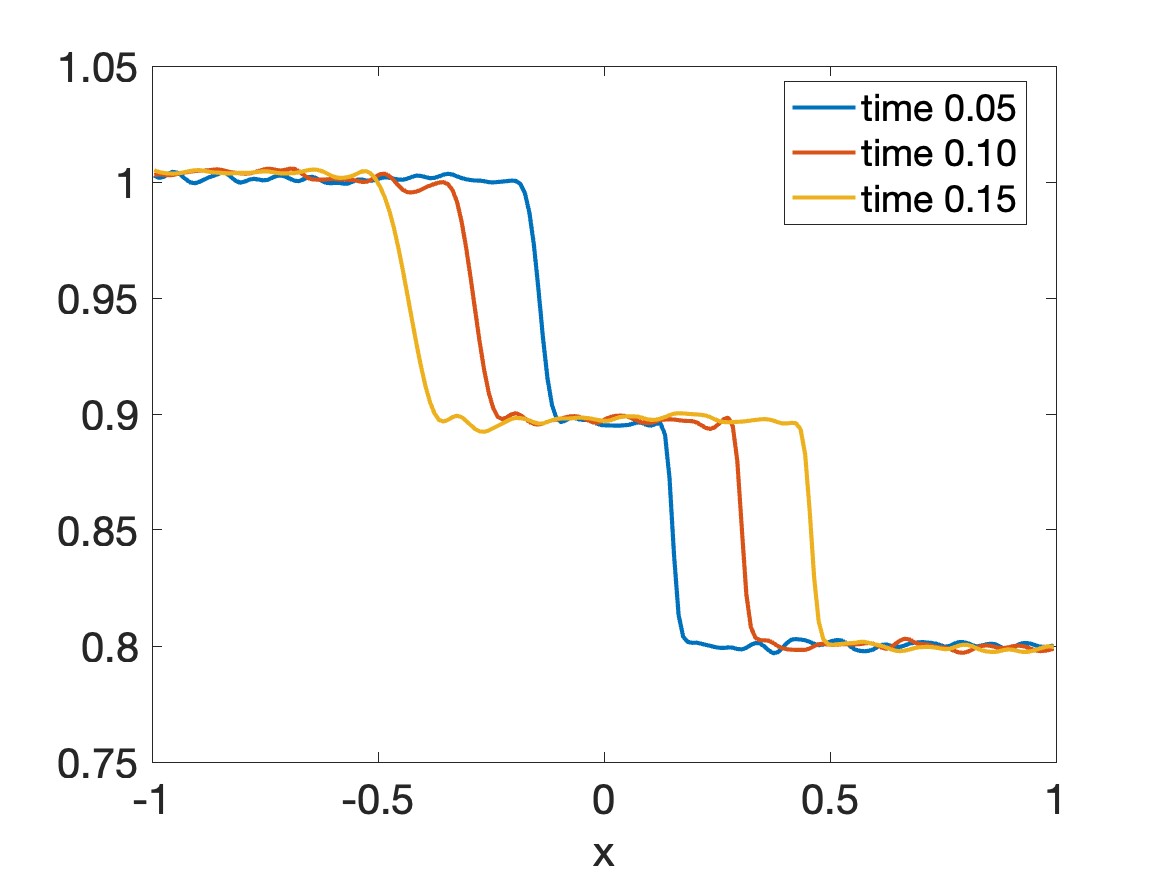}
\includegraphics[width=0.32\textwidth]{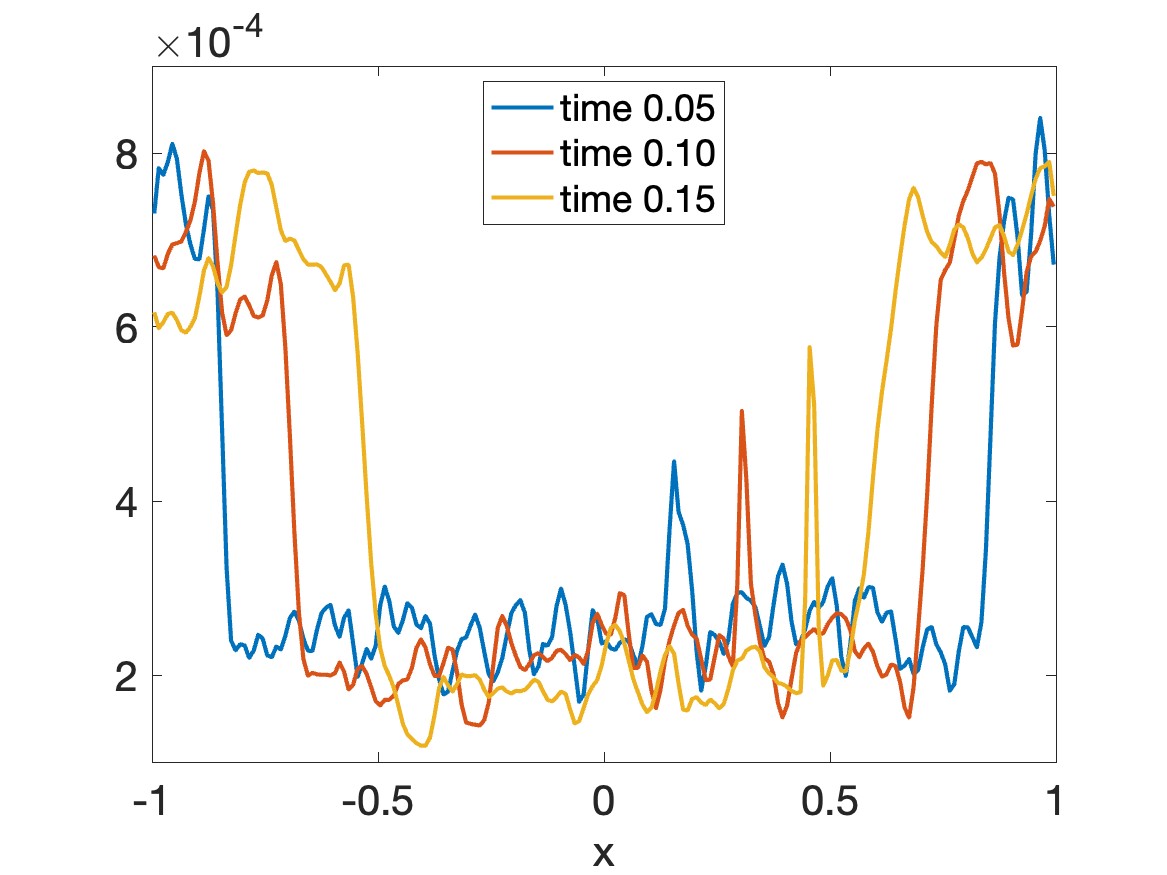}
\includegraphics[width=0.32\textwidth]{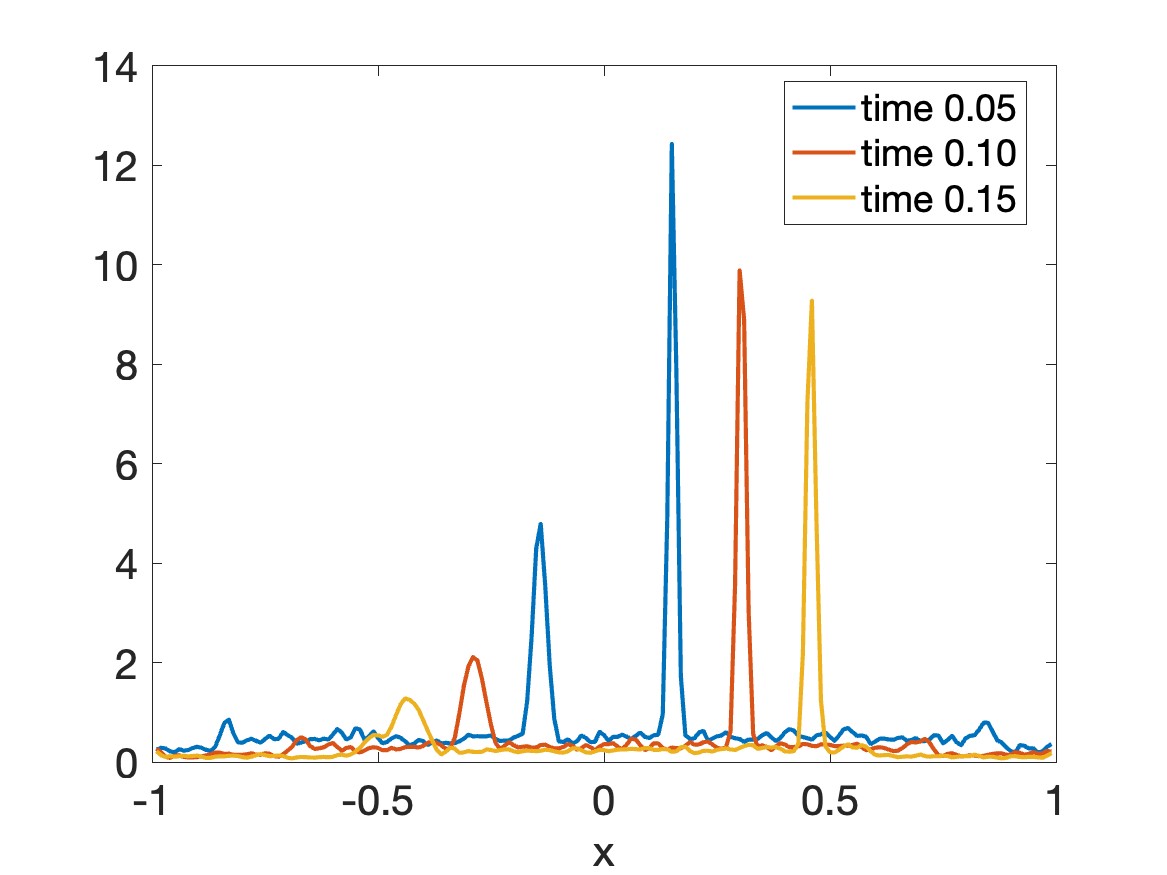}
\caption{The mean solution of $h(x,t)$, $t = 0.05, 0.10, 0.15$ in \eqref{eq:SWE} using fifth order WENO with $\Delta x=10^{-2}$ and $\Delta t=10^{-5}$ for $K = 100$ perturbed initial conditions (left).  The corresponding variance (middle) and gradient second moment (right).}
\label{fig:meanvargrad2}
\end{figure}

Figure \ref{fig:meanvargrad2} helps to provide insight on what information should be used to construct $W$. Figure \ref{fig:meanvargrad2} (left) shows the mean of $K = 100$ numerical solutions for the depth $h(x,t)$ in the shallow water equations \eqref{eq:SWE} for $t = 0.05, 0.10$ and $0.15$. Each numerical solution is computed using fifth order finite difference WENO and is based on perturbed initial conditions in \eqref{eq:initialeasy} with standard deviation $0.05$. Figure \ref{fig:meanvargrad2} (middle) shows the corresponding variance of the solution, while Figure \ref{fig:meanvargrad2} (right) displays the second moment of the gradient of the solution. Observe that while the solution mean exhibits discontinuities and steep gradients at each time instance,  the variance over the whole domain is of $\mathcal O(8\times 10^{-4})$ and as such does not delineate  local structural information.  By contrast,  the gradient second moment not only extracts each discontinuity location, but also apparently distinguishes between shocks and steep gradients, suggesting its potential benefits in determining how heavily {\em each spatial} observation term in \eqref{eq:obj_w} should be penalized. 

\begin{remark}
\label{rem:discontinuityprior}
We note that there are variational data assimilation methods that incorporate structural information into the prior.  For example, the objective function in 4DVar\footnote{The four-dimensional variational data assimilation (4DVar) scheme is a variational data assimilation scheme that optimizes an objective function (similar to the one in Algorithm \ref{alg:ETKF}) based on the given model and observation data \cite{Asch16}.} can be reformulated as an $\ell_2$ regularization \cite{Johnson05}. The method in \cite{Freitag10} uses $\ell_1$ regularization in the analysis step to recover sharp features in the solution. The $\ell_1$ regularization variational approach can also be suitably modified to promote sparsity in the gradient of the state variable \cite{Ebtehaj13, Foufoula-Georgiou14}. A mixed $\ell_1$-$\ell_2$ regularization framework was adopted for the method introduced in \cite{Freitag13, Asadi19}, where the $\ell_1$ norm of the gradient of the analysis vector is added to constrain the  standard $\ell_2$ objective function. By construction, each of these methods use a {\em global} penalty parameter, and therefore does not incorporate {\em local} information into the prior.  Our approach, by contrast, uses the filtering problem for statistical data assimilation.  We furthermore directly incorporate local structural information into the weighting matrix design.  Finally, since our method employs $\ell_2$ minimization, there is a closed-form solution. 
\end{remark}

\subsection{Designing the prior weighting matrix $W$ for densely sampled observations}
\label{sec:full}

In the ideal environment, where the sampled spatial observations are dense, the observation operator $H$ in \eqref{eq:obj_w} is simply the identity matrix $I$. The most straightforward way to construct $W$ in \eqref{eq:obj_w} is then to use the (inflated and localized) ETKF covariance stemming from \eqref{eq:inflation} and \eqref{eq:localization}. We can reasonably choose $\cT = I$ in \eqref{eq:localization}, leading to a diagonal weighting matrix (thereby increasing computational efficiency) given by 
\begin{equation}\label{eq: diag cov inflate}
W=\alpha^2\wh{C}^{D}, \quad\quad \alpha > 1,
\end{equation}
where
\begin{equation*}
\wh{C}^{D}=
\begin{pmatrix}
\wh{C}_1 & 0 & \cdots  & 0 \\
0 & \wh{C}_2 & \cdots  & 0 \\
\vdots & \vdots & \ddots & \vdots \\
0 & 0 & \cdots & \wh{C}_{n} 
\end{pmatrix},
\end{equation*}
with
\begin{equation}
\wh{C}_{i}=\frac{1}{K-1} \sum_{k=1}^{K} (\wh{v}_{i}^{(k)}-\wh{m}_i)^2, \quad i=1, \cdots, n, 
\label{eq: variance}
\end{equation}
being the sample variance of the state variable at each location derived from \eqref{eq: EnKF prior 2}. Here $\wh{v}_{i}^{(k)} = \wh{\bv}^{(k)}(x_i)$ and $\wh{m}_{i}=\wh{\m}(x_i)$, with the subscript denoting the location (spatial) index. In some sense, using \eqref{eq: diag cov inflate} to define $W$ may be viewed as a way to incorporate spatial information into the posterior update, since it describes how well the prior mean approximates the ensembles at each grid point. This construction does {\em not}, however, separately take into account when the underlying solution admits a discontinuity. To this end it is well understood in other applications (e.g.~TV regularization for image recovery) that incorporating such information into the prior typically improves the quality of the numerical solution. Moreover, as was discussed in Section \ref{sec:introduction}, joint information from multiple observations (MMVs) can be used to determine how much to spatially vary the weight of the regularization term in the corresponding objective function. Drawing an analogy between the weighted $\ell_1$ joint sparsity image recovery algorithm in \cite{Adcock19,GelbScarnati2019}, and the objective function in \eqref{eq:obj_w}, where the ensemble corresponds to the MMVs in the image recovery problem,  we are inspired to construct $W$ using joint information from the ensemble data so that the optimized solution is one that is able to  distinguish discontinuous from smooth regions. Consistent to what is observed in Figure \ref{fig:meanvargrad2}, we will use information from the {\em spatial gradient} of the state variable, which is sparse for a piecewise constant solution profile, and is therefore better suited for designing $W$ than the state variable is itself.  We now provide details for this construction. 

In our problem we consider the grid to be uniform with mesh size $\Delta x$. For simplicity we compute the gradient using first order finite differencing,
\begin{equation}
\label{eq:gradensembles}
d\wh{v}^{(k)}_{i+\frac{1}{2}}=\frac{\wh{v}^{(k)}_{i+1}-\wh{v}^{(k)}_{i}}{\Delta x},
\end{equation}
and note that neither the uniformity of the grid nor the approximation of the gradient are intrinsic to our approach.  The second moment of the gradient at each cell center is then 
\begin{equation}
\label{eq:secmomentcenter}
\wh{S}_{i+\frac{1}{2}}=\frac{1}{K} \sum_{k=1}^{K} \left( d\wh{v}^{(k)}_{i+\frac{1}{2}} \right)^2.
\end{equation}
We use \eqref{eq:secmomentcenter} to approximate the second moment of the gradient at each grid point $\{x_i\}_{i = 1}^n$ as 
\begin{equation}\label{eq: second moment}
\wh{S}_i = \begin{cases}
\ds \frac{1}{2}\wh{S}_{\frac{3}{2}},\qquad i=1,\\[2ex]
\ds \frac{1}{2}(\wh{S}_{i-\frac{1}{2}}+\wh{S}_{i+\frac{1}{2}}),\qquad i=2, \cdots, n-1\\[2ex]
\ds \frac{1}{2}\wh{S}_{n-\frac{1}{2}}, \qquad i=n,
\end{cases}
\end{equation}
from which we construct the weighting matrix 
\begin{equation}\label{eq: diag gsm}
W = \beta \wh{S}^{D},\quad \beta > 0,
\end{equation}
with 
\begin{equation*}
\wh{S}^{D}= 
\begin{pmatrix}
\wh{S}_1 & 0 & \cdots  & 0 \\
0 & \wh{S}_2 & \cdots  & 0 \\
\vdots & \vdots & \ddots & \vdots \\
0 & 0 & \cdots & \wh{S}_{n} 
\end{pmatrix}.
\end{equation*}

\subsection{Designing the prior weighting matrix $W$ for sparsely sampled observations}
\label{sec:sparse}

As already noted, using a diagonal weighting matrix $W$ is suitable when observations are densely sampled.  We now seek a more flexible design for $W$ to accommodate a sparse observation environment.
Based on $\wh{C}_i$ in \eqref{eq: variance}, the full prior covariance matrix $\wh{C}$ has entries
\begin{equation}
\label{eq:full_covariance}
\wh{C}_{i,i'}=\sqrt{\wh{C}_{i}}\,\wh{r}_{i,i'}\sqrt{\wh{C}_{i'}}, \quad i,i' = 1,\dots,n,
\end{equation}
where $\wh{r}_{i,i'}$ given by 
\begin{equation}\label{eq: correlation coefficient}
\wh{r}_{i,i'} = \frac{\ds \frac{1}{K-1} \sum_{k=1}^{K} (\wh{v}_{i}^{(k)}-\wh{m}_i)(\wh{v}_{i'}^{(k)}-\wh{m}_{i'})}{\ds \sqrt{\frac{1}{K-1} \sum_{k=1}^{K} (\wh{v}_{i}^{(k)}-\wh{m}_i)^2} \sqrt{\frac{1}{K-1} \sum_{k=1}^{K} (\wh{v}_{i'}^{(k)}-\wh{m}_{i'})^2}},
\end{equation}
are the components of the correlation coefficient matrix $\wh{R}$.
The idea is then to generalize  \eqref{eq: diag gsm} by making use of the fact that the entries in the full sample covariance matrix $\wh{C}$ in \eqref{eq:full_covariance} are related to the sample variance $\{\wh{C}_i\}_{i = 1}^n$ in \eqref{eq: variance} through $\wh{R}$.
That is, $\wh{C}$, which is equivalently formed from either \eqref{eq:full_covariance} or \eqref{eq: EnKF prior 2} (where the subscript $j+1$ is a {\em temporal} index, not to be confused with the {\em spatial} indices used in this section), can be directly determined from \eqref{eq: variance}. 

Based on the comparison of  \eqref{eq: diag cov inflate} to \eqref{eq: diag gsm}, we are similarly motivated to now construct a full gradient second moment matrix that incorporates the gradient second moment components given in \eqref{eq: second moment} with the correlation coefficients in \eqref{eq: correlation coefficient}. That is, we directly substitute $\wh{S}_i$ for $\wh{C}_i$ in \eqref{eq:full_covariance} to obtain
$\wh{S}$ 
\begin{equation}\label{eq: full second gradient}
\wh{S}_{i,i'}=\sqrt{\wh{S}_{i}}\,\wh{r}_{i,i'}\sqrt{\wh{S}_{i'}}.
\end{equation}
As before, we construct the weighting matrix as 
\begin{equation}\label{eq: new weighting matrix}
W = \beta \wh{S}\odot \cT,
\end{equation}
with  $\beta > 0$ and $\cT$ the user prescribed binary Toeplitz localization matrix as is used in \eqref{eq:localization}.

\begin{remark}
Based on the observation that \eqref{eq: ETKF C} is equivalent to \eqref{eq:full_covariance} for $\wh{X} = (\wh{C}^D)^{\frac{1}{2}} \wt{X}$ when $\wt{X}$ is given by $\wh{R}=\wt{X} \wt{X}^T$, rather than using \eqref{eq: full second gradient} in the numerical implementation, we simply compute  $\wh{S} = \left((\wh{S}^D)^{\frac{1}{2}} \wt{X} \right) \left((\wh{S}^D)^{\frac{1}{2}} \wt{X} \right)^T$. 
\end{remark}

\subsection{Structural refinement of the prior by clustering}
\label{subsec: clustering}

Intuitively speaking, three types of information can be used to update the state variable at any particular location: (1) the observation data at that location; (2) the prior information at that location; and (3) relevant information coming from other (neighboring) locations. As discussed in Section \ref{sec:full}, we can effectively rebalance these first two types of information using the gradient information via \eqref{eq: diag gsm}, while incorporating neighborhood information was described in Section \ref{sec:sparse} through use of \eqref{eq: full second gradient}.  Utilizing such correlating information becomes increasingly important  in cases where we lack observation data at that particular location.   Moreover, when the underlying state variable has a discontinuity profile, gradient information provides a more appropriate weighting matrix, and in particular promotes solutions that ``trust'' prior information  more. Intrinsically, then, there is an underlying assumption that the prediction method used to simulate the dynamic model is well-suited to solve PDEs that contain discontinuities in their solutions. We will return to this discussion in Section \ref{sec:model}.

It is also possible to update the correlation matrix $\wh{R}$ with the entries given by \cref{eq: correlation coefficient} that may help to refine the {\em local} information regarding discontinuities. In particular, the state variable at locations on either side of a discontinuity should be less correlated. We now propose a clustering approach to extract this information and improve our correlation update. We note that our clustering approach is similar to localization (Section \ref{sec:inflatelocal}), in the sense that both aim to reduce unwanted correlations informed by physical properties.

We begin by using the gradient information of the prior mean \eqref{eq: EnKF prior 1} to find the jump discontinuity locations, specifically by computing
\begin{equation}
d\wh{m}_{i+\frac{1}{2}}=\frac{\wh{m}_{i+1}-\wh{m}_{i}}{\Delta x}, 
\label{eq:gradmean}
\end{equation}
from which the local extrema determine the set of discontinuity locations. While there are several ways to do this, assumptions regarding both the underlying problem and the available data are always required. For example, one could impose constraints according to distance between jump locations with respect to grid point resolution. Consistent with our numerical example in Section \ref{sec:numerical}, for ease of presentation and implementation here we assume that there is only one true discontinuity. Thus we only seek to find the global extrema in \eqref{eq:gradmean}, and define the corresponding discontinuity location as 
\begin{equation}
\label{eq:discpoints}
\xi = \underset{i}{\text{argmax}} \{ |d\wh{m}_{i+\frac{1}{2}}| \}. 
\end{equation}
We emphasize that neither the number of discontinuities nor the way we determine their locations inherently limits the use of our method, which can indeed consider an unknown number of discontinuities along with more sophisticated (e.g.~higher order) ways of determining their locations. There would be additional considerations for computational complexity, however.
Finally, we note that we choose to use \eqref{eq:gradmean} rather than \eqref{eq:gradensembles}, which is more likely affected by noise. 

Once $\xi$ in \eqref{eq:discpoints} is determined, we can define the discontinuity region
\begin{equation}\label{eq: dist}
R^D = \{ i: |i-\xi| \leq dist \}
\end{equation}
where $dist \in \mathbb{Z}^+$ is predetermined to take into account both gradient steepness and the resolution of the grid. The domain is then decomposed into the discontinuity region and two smooth regions
$$ R^S_1 = \{ i: i < i' \quad \forall i' \in R^D  \} \qquad \text{and} \qquad  R^S_2 = \{ i: i > i' \quad \forall i' \in R^D  \}. $$
The entries of the correlation coefficient matrix in \eqref{eq: correlation coefficient}, are modified according to each entry's proximity to a discontinuity region as
\begin{equation}
\label{eq:mod_coeffs}
\wt{r}_{i,i'} = 
\begin{cases}
\wh{r}_{i,i'} &\quad \text{if } i=i' \text{ or } i,i' \in R^S_\kappa, \, \kappa=1,2, \\
0 &\quad \text{otherwise}, \\
\end{cases}    
\end{equation}
leading to the modified correlation coefficient matrix $\wt{R}$.
By construction, \eqref{eq:mod_coeffs} imposes local information so that only values corresponding to spatial locations occurring in the {\em same} smooth region can be correlated, or equivalently, values corresponding to spatial locations on opposite sides of a discontinuity region will not be correlated. Moreover, the solution within a discontinuous region $R^D$ will not utilize information from neighboring locations.
Replacing the correlation coefficients in \eqref{eq: full second gradient} with \eqref{eq:mod_coeffs} we obtain the block diagonal matrix (thereby also improving computational efficiency)
\begin{equation}\label{eq: full second gradient 2}
\wt{S}_{i,i'}=\sqrt{\wh{S}_{i}}\,\wt{r}_{i,i'}\sqrt{\wh{S}_{i'}},
\end{equation}
from which we compute 
\begin{equation}\label{eq: new weighting matrix 2}
W = \beta \wt{S}\odot \cT,
\end{equation}
as in \eqref{eq: new weighting matrix} for some parameter $\beta$ and localization matrix $\cT$.

\subsection{Numerical implementation}

Combining the results from each previous section, we now recap the modified ETKF method using the new second moment weighting matrix, $W = \beta \wt{S} \odot \cT$, where $\wt{S}$ is given by \eqref{eq: full second gradient 2}, with  $\wt{S} = \wh{S}$ in \eqref{eq: full second gradient} when clustering is not included.  We note again that both the tuning parameter $\beta > 0$ and the  localization  matrix $\cT$ are user prescribed.  In our numerical experiments, we choose $\cT= I$ when observations are densely sampled (Section \ref{sec:full}) and  to be tridiagonal with nonzero entries of $1$ when the observations are sparse (Section \ref{sec:sparse}).  At each temporal increment $t = t_j$, the prior ensembles $\wh{\bv}^{(k)}$ are obtained through the dynamic model  \eqref{eq: EnKF predict 1}. The sample prior mean $\wh{\m}$ is computed as \eqref{eq: EnKF prior 1}. The prior covariance matrix 
$$\wh{C} = \wh{X}\wh{X}^T$$
is generated using the non-inflated centered ensemble
\begin{equation}
\label{eq:noninflated}
\wh{X} = \frac{1}{\sqrt{K-1}} \left[ \wh{\bv}^{(1)}-\wh{\m}, \cdots, \wh{\bv}^{(K)}-\wh{\m} \right].
\end{equation}
We then construct $T$ and $X$ respectively according to \eqref{eq: ETKF T} and \eqref{eq: ETKF X}, followed by $\wt{S}$ as \eqref{eq: full second gradient 2} or \eqref{eq: full second gradient}, depending on whether or not clustering is used. To determine the posterior mean in the analysis step, we now solve \eqref{eq:obj_w} as
\begin{equation}
\label{eq:ETKF second moment}
\m = \underset{\bv}{\text{arg min}} \left( \frac{1}{2}\vert \by-H \bv \vert_{\Gamma}^2 + \frac{1}{2} \vert \bv-\wh{\m} \vert_{\beta \wt{S}\odot \cT}^2\right),
\end{equation}
where we have written $\beta \wt{S}\odot \cT$ in place of $W$ since $\beta$ and $\cT$ are separately prescribed.
This is followed by updates on the ensembles according to the transformation
\begin{gather}
\bv^{(k)}=\m+ \bzeta^{(k)},
\end{gather}
where $\bzeta^{(k)}$ is defined by \eqref{eq: ETKF zeta}.  Observe that \eqref{eq:ETKF second moment} can be analytically solved for $\m$ via
\begin{equation*}
\left(H^T \Gamma^{-1} H + \left( \beta \wt{S}\odot \cT \right)^{-1} \right) \m = H^T \Gamma^{-1}\by + \left( \beta \wt{S}\odot \cT \right)^{-1} \wh{\m}.    
\end{equation*}
Alternatively, to avoid direct computation of $(\beta \wt{S}\odot \cT)^{-1}$ we use the equivalent Sherman-Morrison-Woodbury formula yielding 
\begin{subequations}\label{eq:SMWformula}
\begin{equation}
    \m=(I-KH)\wh{\m}+ K \by, \label{eq: aETKF post 3}
    \end{equation}
\begin{equation}
K=\beta \wt{S}\odot \cT H^{T}S^{-1}, \quad S=H \left( \beta \wt{S}\odot \cT \right) H^{T} + \Gamma.\label{eq: aETKF post 4}
\end{equation}
\end{subequations}

Algorithm \ref{alg:mETKF} describes the modified ETKF framework using the new weighting matrix $W = \beta \wt{S}\odot \cT$ with localization and clustering techniques.

\begin{algorithm}[ht!]
\caption{modified ETKF method with localization and clustering techniques}
\begin{algorithmic}
\State Initialize ensembles $\bv^{(k)}_0$, $k = 1,\dots,K$. Choose parameter $\beta$ and localization matrix $\cT$ for $W$ in \eqref{eq: new weighting matrix 2}.
\For{$j=1$ to $J-1$} 
\State Compute prior ensembles $\ds \wh{\bv}^{(k)}_{j+1}=\bPsi(\bv^{(k)}_j)$ for $k=1, \cdots, K.$
\If{data observations are available at time instance $j+1$}
\State Compute prior mean $\ds \wh{\m}_{j+1}=\frac{1}{K}\sum_{k=1}^{K}\,\wh{\bv}^{(k)}_{j+1}$.
\State Compute non-inflated centered ensemble $\ds \wh{X}_{j+1} = \frac{1}{\sqrt{K-1}} \left[ \wh{\bv}^{(1)}_{j+1}-\wh{\m}_{j+1}, \cdots, \wh{\bv}^{(K)}_{j+1}-\wh{\m}_{j+1} \right]$. 
\State Compute the prior covariance matrix $\ds \wh{C}_{j+1}=\wh{X}_{j+1}\,\wh{X}_{j+1}^T$.
\State Compute the transformation matrix $T_{j+1}$ and $X_{j+1}$ according to \eqref{eq: ETKF T} and \eqref{eq: ETKF X} respectively.
\State Compute the full gradient second moment matrix $\wt{S}_{j+1}$ according to \eqref{eq: full second gradient 2}.
\State Solve the minimization problem for the posterior mean $\m_{j+1}$ as
\begin{gather*}
\m_{j+1}=(I-KH)\wh{\m}_{j+1}+ K \by_{j+1}, 
\end{gather*}
where
\[K=\beta \wt{S}_{j+1} \odot \cT H^{T}S^{-1},\quad\quad S=H \left( \beta \wh{S}_{j+1}\odot \cT \right) H^{T} + \Gamma.\]
\State Update the posterior ensembles according to the transformation \eqref{eq:zeta}
where $\bzeta^{(k)}_{j+1}$ is defined by \eqref{eq: ETKF zeta}.
\Else{ $\ds \bv^{(k)}_{j+1}=\wh{\bv}^{(k)}_{j+1}$ for $k=1, \cdots, K.$}
\EndIf
\EndFor
\end{algorithmic}
\label{alg:mETKF}
\end{algorithm}

\begin{remark}
\label{rem:inflatingvsparametertuning}
In contrast to traditional ETKF,  where inflation is directly applied to the centered ensembles via \eqref{eq:inflation}, Algorithm \ref{alg:mETKF} uses the non-inflated ensemble \eqref{eq:noninflated} to compute the covariance matrix which is in turn employed to construct the correlation coefficient matrix with entries given by \eqref{eq: correlation coefficient}. A tuning parameter,  here $\beta$, is still needed for the weighting matrix construction, however. 
\end{remark}

\section{Numerical examples}
\label{sec:numerical}
We use the one-dimensional dam break problem governed by shallow water equations as our model problem to demonstrate the performance of our new method.  The dam break problem yields a solution that has both shock and contact discontinuities, making it challenging for traditional statistical data assimilation techniques.  By contrast, our approach is designed especially for this class of problems. The analytical solution provided in \cite{Stoker} conveniently allows us to conduct some error analysis. While for simplicity this investigation only considers the 1D case, our method is not inherently limited, and 2D problems will be discussed in future work.

\subsection{1D dam break problem}
\label{sec:model}

Assuming a flat bottom, the 1D shallow water equations in conservation form for spatial variable $x$ and temporal variable $t$ is given by
\begin{subequations}
\label{eq:SWE}
\begin{equation}
h_t+(hu)_x=0, \label{eq: SWE1} 
\end{equation}
\begin{equation}
(hu)_t+(hu^2+\frac{1}{2}gh^2)_x=0, \label{eq: SWE2}
\end{equation}
\end{subequations}
where $h=h(x,t)$ is the depth of the water, $u=u(x,t)$ is the (depth-averaged) fluid velocity, and $g = 9.81$ is the acceleration constant due to gravity. 

For the 1D dam break problem, the fluid is initially at rest on both sides of an infinitely thin dam located at $x=0$. A dam break is simulated by sudden removal of the dam wall at time $t = 0$. Stoker \cite{Stoker} provides an analytical solution for this problem with initial conditions given by
\begin{gather}
\label{eq:initialeasy}
h(x,0)=
\begin{cases}
h_0, \quad x<0 \\
h_1, \quad x>0
\end{cases} \qan
u(x,0)=0.
\end{gather}
The solution consists of a shock front propagating to the right and a rarefaction wave propagating to the left. This analytical solution will serve as our underlying truth for the test problems in Section \ref{subsection:numerical full} and Section \ref{subsection:numerical sparse} and will be referred to as the ``Stoker solution". 
 
For the experiments conducted in Section \ref{subsection:numerical full} and Section \ref{subsection:numerical sparse}  we let $h_0=1$, $h_1=0.8$,  while Section \ref{subsection:numerical oscillate} we consider oscillatory initial conditions.  The spatial domain is always $[-1,1]$.  To set up our experiments, we first numerically simulate the coupled system in \eqref{eq:SWE} and then use the corresponding solution for the velocity $u(x,t)$ as the synthetic data for the data assimilation on recovering the water depth, $h(x,t)$.  

It is, of course, important to choose an appropriate numerical PDE solver. Since our goal is to analyze the impact of our new weighting matrix, here we employ a high-order method to minimize the error coming from the PDE simulation.  However, it is well known that for solutions containing discontinuities, traditional high order finite difference schemes typically yield spurious oscillations near discontinuities, which may pollute smooth regions and even lead to instability, causing numerical blowup \cite{Shu20}. The WENO method \cite{Liu94} is designed to achieve higher-order accuracy in smooth regions while sharply resolving discontinuities in an essentially non-oscillatory fashion. We therefore use fifth order finite difference WENO spatial discretization along with Lax-Friedrichs flux splitting in our numerical experiments.\footnote{For the Lax-Friedrichs flux splitting parameter, we simply choose the default value $\lambda =\max(|u|+\sqrt{gh})$.} Finally, we emphasize that other numerical PDE solvers can be straightforwardly accommodated, and the choice of method should be based on the usual considerations (e.g.~domain regularity, grid resolution, computational efficiency).  

To close the system we apply homogeneous Dirichlet boundary conditions ($u(\pm 1,t) =0$) at both ends of the domain, noting that since the boundaries are artificial, the numerical boundary conditions are prescribed simply to ensure that there is no reflection. We also use the total variation diminishing third-order Runge-Kutta (TVDRK3) time integration scheme \cite{Liu94} with a time step that satisfies the CFL condition. Finally, we remark that no attempt was made to optimize parameters with regard to either increased resolution or efficiency.




To simplify our presentation, our numerical experiments only uses \eqref{eq: SWE1}, the depth of the water $h$ advected by the transport equation, for the dynamic model \eqref{model: dynamic}. We assume that the fluid velocity $u(x,t)$ is known, with its value at any point in the domain taken from the numerical solution of the coupled problem. The domain $[-1,1]$ is uniformly discretized over $n = 1001$ points, with  $x_i = -1 + (i-1)\Delta x$, $i = 1,\dots,n$, yielding  $\Delta x = \frac{2}{n-1} = 2 \times 10^{-3}$.  We use $\Delta t = CFL\# \Delta x$ to integrate in time, with $CFL\# = 0.1$ to ensure the CFL condition is satisfied  (so that $\Delta t = 2 \times 10^{-4}$). The observation data are created from the Stoker solution with $\Gamma = \gamma^2 I$, $\gamma = 0.01$ in \eqref{assump: noise Gaussian}, and are observable every $\Delta t_{obs} = 5 \Delta t = 10^{-3}$. 
We further consider two observational data environments: (1) The {\em densely sampled} case, where observations are made at every point in the spatial domain;  and (2) the {\em sparsely sampled} case, where observations are available at every other point. Finally, the $K=100$ ensembles
in Algorithm \ref{alg:mETKF} are obtained from the data assimilation solution initialized from either \eqref{eq:initialeasy} or \eqref{eq: init}, and in both cases contain additive i.i.d.~noise with mean $0$ and standard deviation $0.1$. 

\begin{remark}
\label{rem:ensemblesize}
In practical applications, the choice of ensemble size depends on dimension reduction and computational feasibility. As our goal here is to compare the performance of methods in Algorithm \ref{alg:ETKF} and Algorithm \ref{alg:mETKF}, we set $K = 100$ to avoid additional fine-tuning.  Supplementary numerical experiments (not reported here) confirmed that results did not fundamentally change for varying values of $K$.
\end{remark}

\begin{figure}[h!]
\centering
\includegraphics[width=0.32\textwidth]{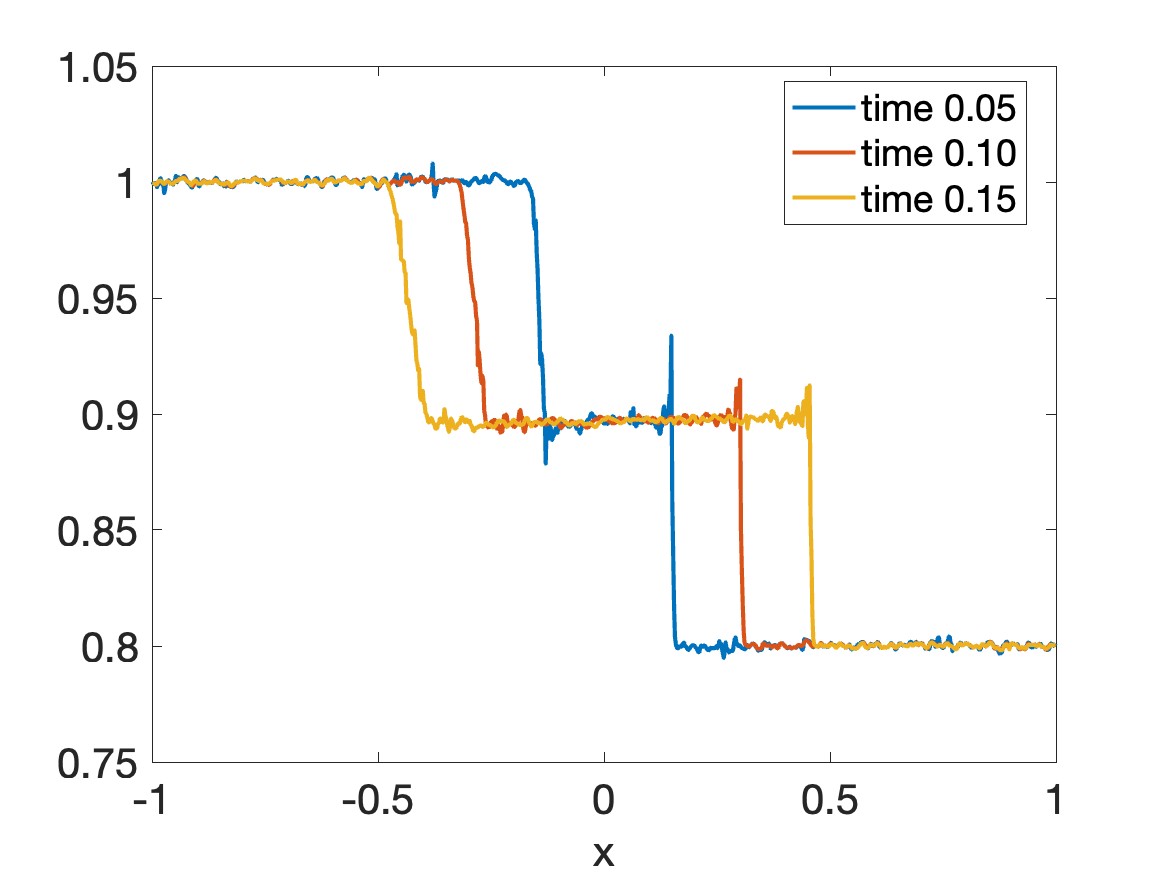}
\includegraphics[width=0.32\textwidth]{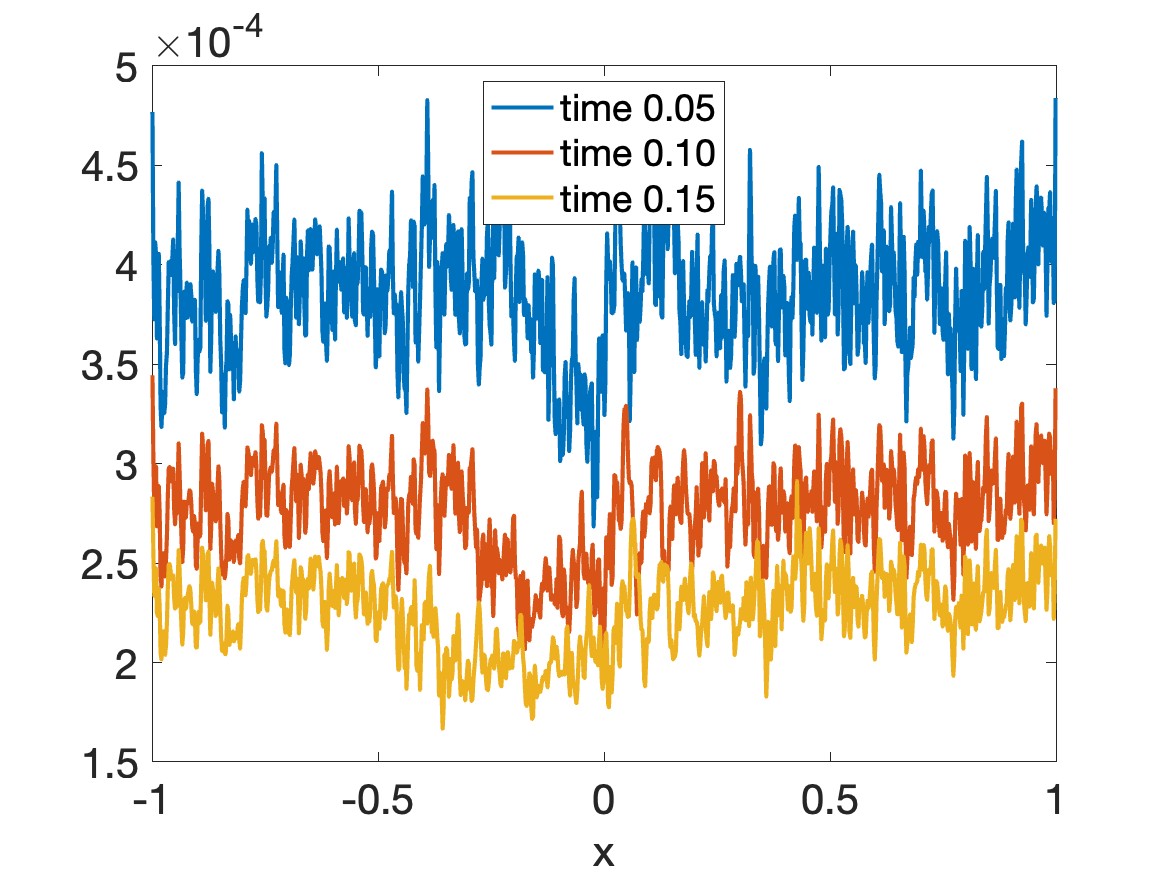}
\includegraphics[width=0.32\textwidth]{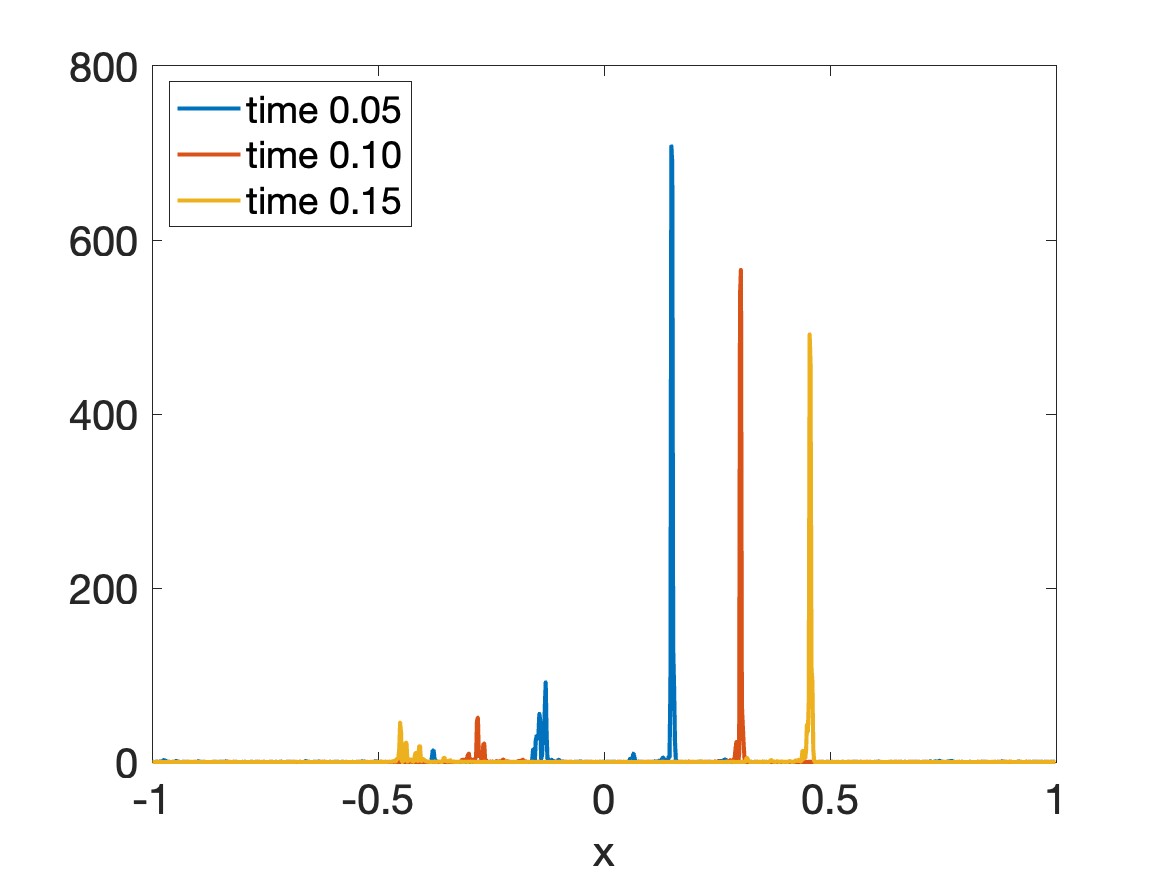}
\caption{The mean solution of the prior ensembles of $h(x,t)$, $t = 0.05, 0.10, 0.15$ as determined by Algorithm \ref{alg:mETKF}.  The corresponding variance (middle) and gradient second moment (right).}
\label{fig:meanvargrad1}
\end{figure}

Similar to the discussion surrounding Figure \ref{fig:meanvargrad2}, Figure \ref{fig:meanvargrad1} demonstrates the benefits of using the gradient second moment of the prior ensembles to generate the weighting matrix in  Algorithm \ref{alg:mETKF}. Specifically, the figure shows the mean, variance and gradient second moment of the prior ensembles for the densely sampled observation case described in Section \ref{sec:full} at time instances  $0.05$, $0.10$ and $0.15$, where the prior ensembles are obtained via Algorithm \ref{alg:mETKF}. Once again we see that although the solution clearly exhibits discontinuities at each time instance,  the variance is not able to delineate any local structural information.  By contrast,  the gradient second moment not only extracts each discontinuity location, but also apparently distinguishes between shocks and steep gradients, suggesting its potential benefits for reweighting the objective function \eqref{eq:obj_w} in order to effectively capture such information  in the solution. More details  follow in Sections \ref{subsection:numerical full} and \ref{subsection:numerical sparse}.

We analyze our results using the pointwise error evaluated as
\begin{equation} 
err(x_i,t) = \vert h(x_i,t) - h_i(t) \vert,\quad i = 1,\dots,n,
\label{eq:error}
\end{equation}
where $h(x_i,t)$ is the Stoker solution at time $T$ for each $x_i$ and $h_i(t)$ is the posterior mean solution obtained either from Algorithm \ref{alg:ETKF} or Algorithm \ref{alg:mETKF},
along with the relative error given by
\begin{equation}
\label{eq:err}
\ds e(t_j) = \frac{\ds \sum_{i = 1}^n |h(x_i,t_j)- h_i(t_j)|}{\ds \sum_{i = 1}^n |h(x_i,t_j)|}
\end{equation}
for each $t_j = j\Delta t$, $j  = 1, \cdots, J$, with $\Delta t = \frac{T}{J}$. 
 
\subsection{Densely sampled observations}\label{subsection:numerical full}

We now compare results obtained using Algorithm \ref{alg:mETKF} with localization matrix $\cT = I$ to those obtained via Algorithm \ref{alg:ETKF} for the ideal densely sampled observation case. 
The inflation parameter $\alpha = 1.5$ comes from hand-tuning to produce a ``good'' root mean square error (RMSE). To determine $\beta$ for the gradient second moment weighting matrix $W=\beta \wh{S}^D$ in \eqref{eq: diag gsm}, we follow the general principle that $\beta$ should be chosen so that $W$ has comparable magnitude to that of the prior covariance matrix used in Algorithm \ref{alg:ETKF}. We therefore select $\beta$ so that $\max_{i,i'} (W_{i,i'}) = 0.003$, $1 \le i,i' \le n$.  As the fine tuning of the parameters is not the focus of our investigation, we do not attempt to optimize them further. 

\begin{figure}[h!]
\centering
\includegraphics[width=0.32\textwidth]{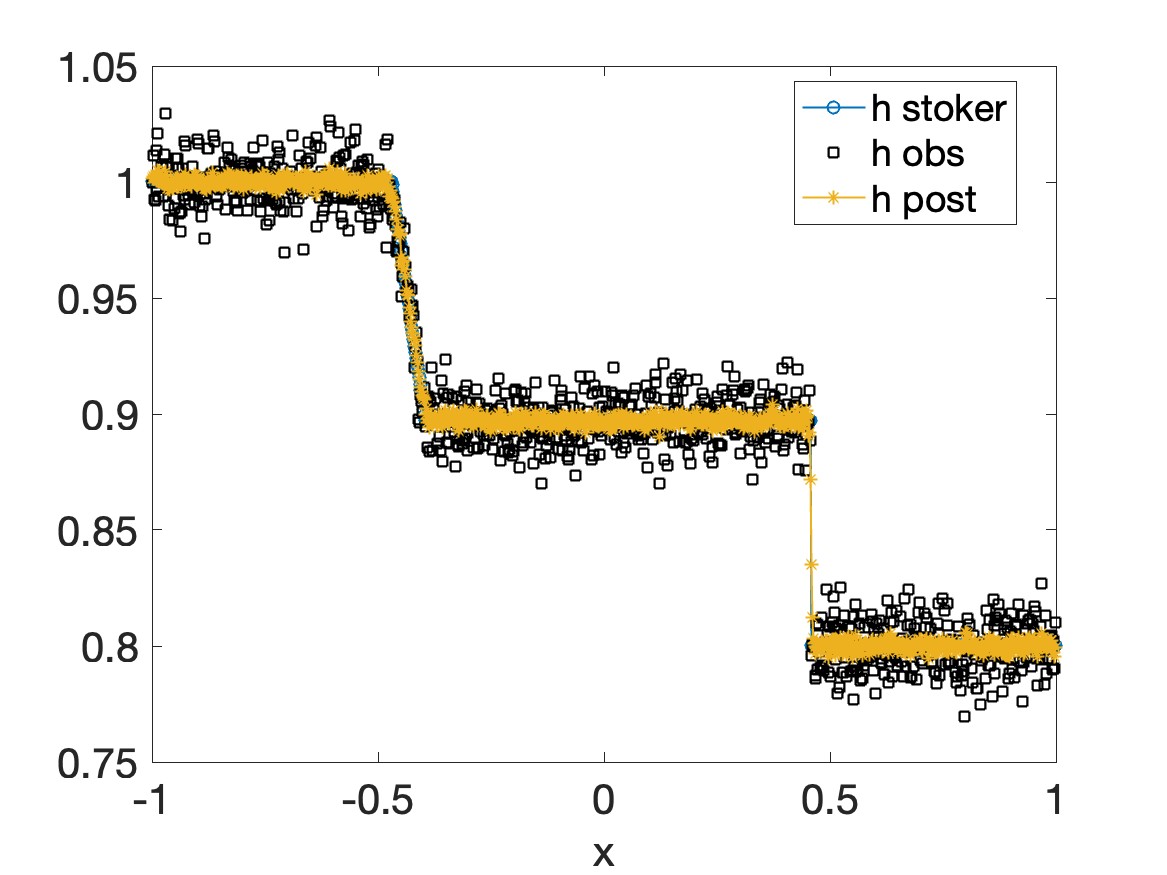}
\includegraphics[width=0.32\textwidth]{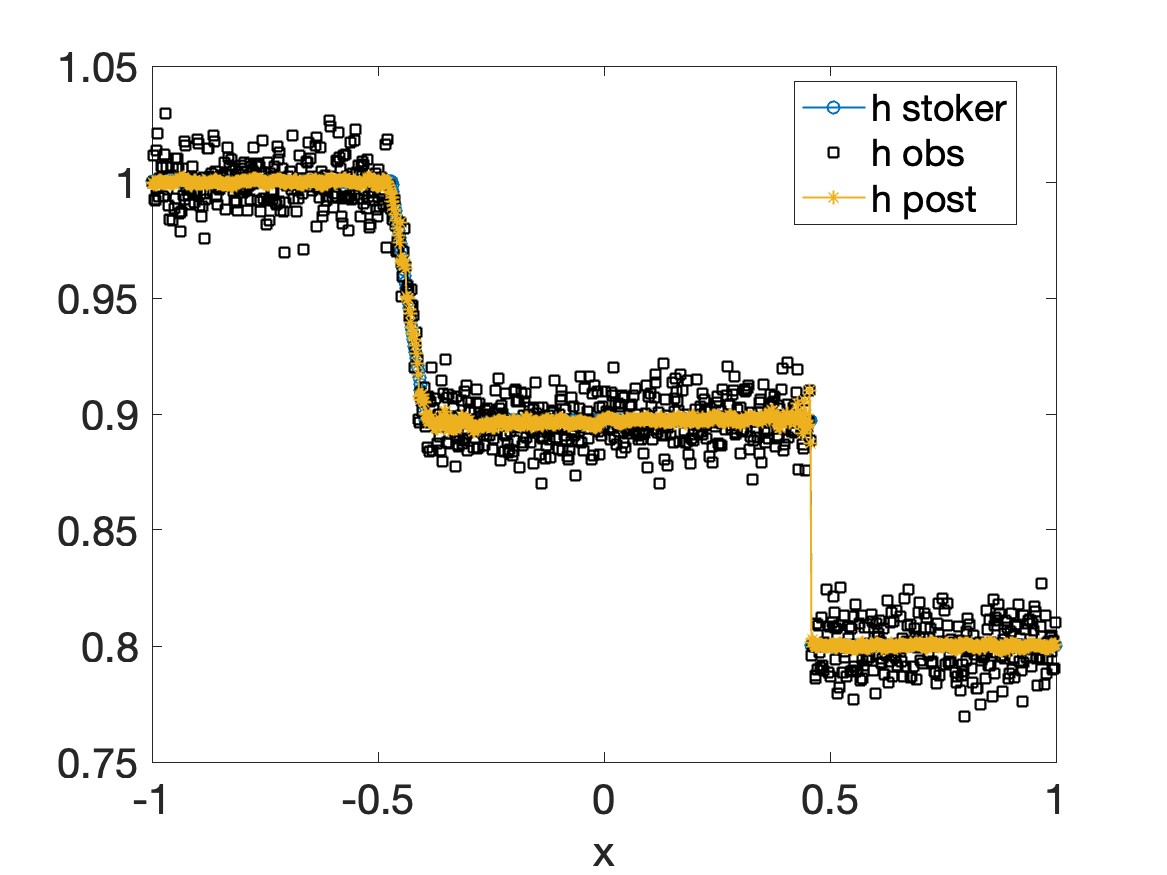}\\
\includegraphics[width=0.32\textwidth]{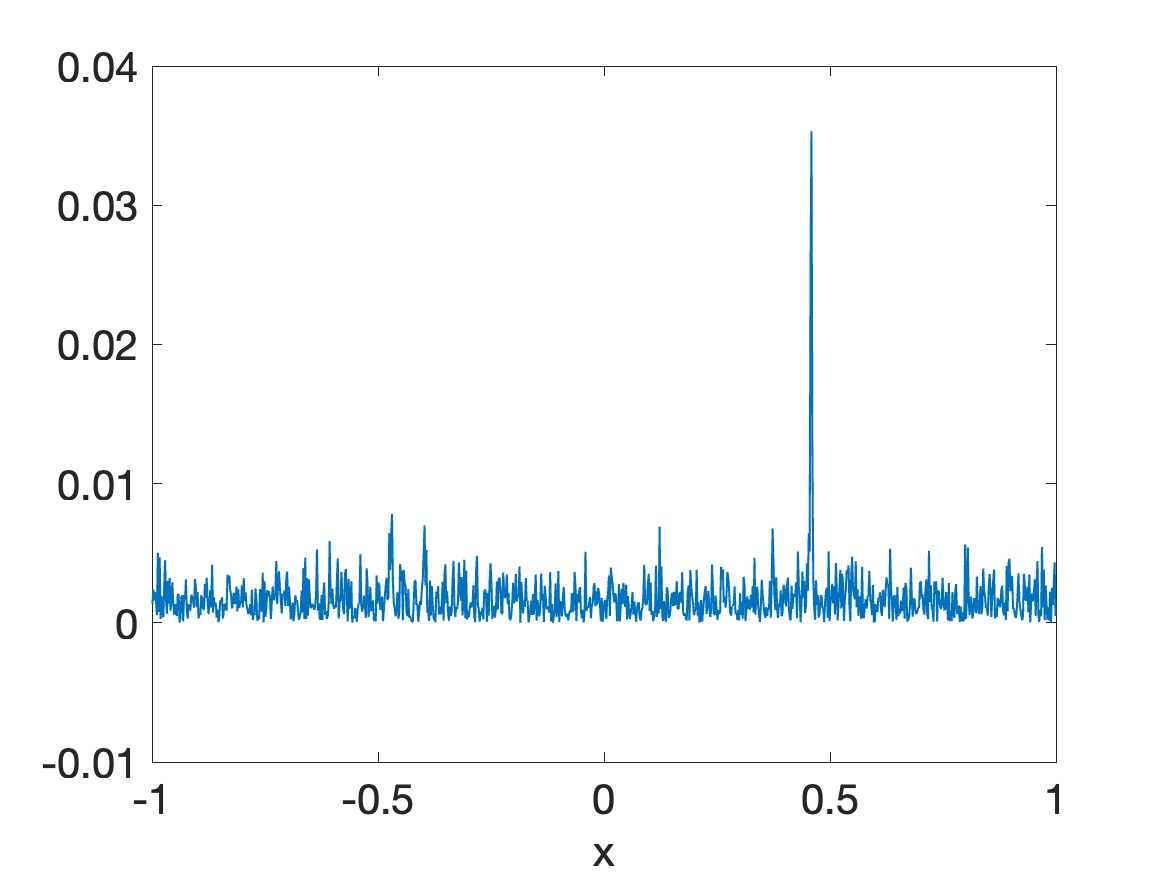}
\includegraphics[width=0.32\textwidth]{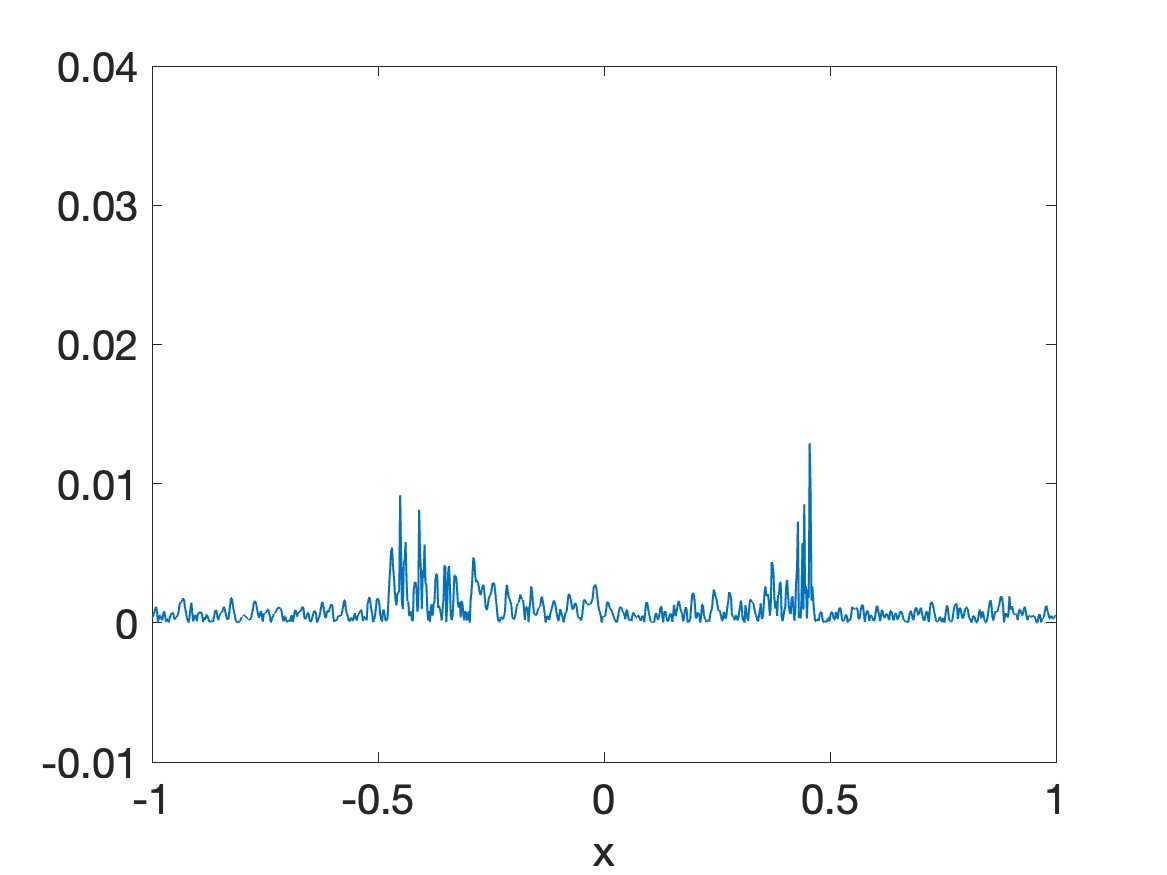}
\caption{(top) Numerical solutions for the posterior mean of depth $h(x,0.15)$ as obtained by (left) Algorithm \ref{alg:ETKF} and (right) Algorithm \ref{alg:mETKF}.  (bottom)  Corresponding pointwise posterior error defined by \eqref{eq:error}. }
\label{Fig: full sol}
\end{figure}

Figure \ref{Fig: full sol}(top) shows the numerical solutions of depth $h(x,0.15)$, along with the Stoker solution and the observational data,  as obtained by  (left) Algorithm \ref{alg:ETKF} and (right)  Algorithm \ref{alg:mETKF}, each using the parameter values as already prescribed. Figure \ref{Fig: full sol}(bottom) displays the corresponding pointwise  posterior errors using \eqref{eq:error}. Observe that employing the gradient second moment weighting matrix not only provides a less oscillatory solution in the smooth regions, but it also generates a more highly resolved shock.  Comparing the magnitudes of the error in Figure \ref{Fig: full sol}(bottom) further verifies these observations.

\begin{figure}[h!]
\centering
\includegraphics[width=0.32\textwidth]
{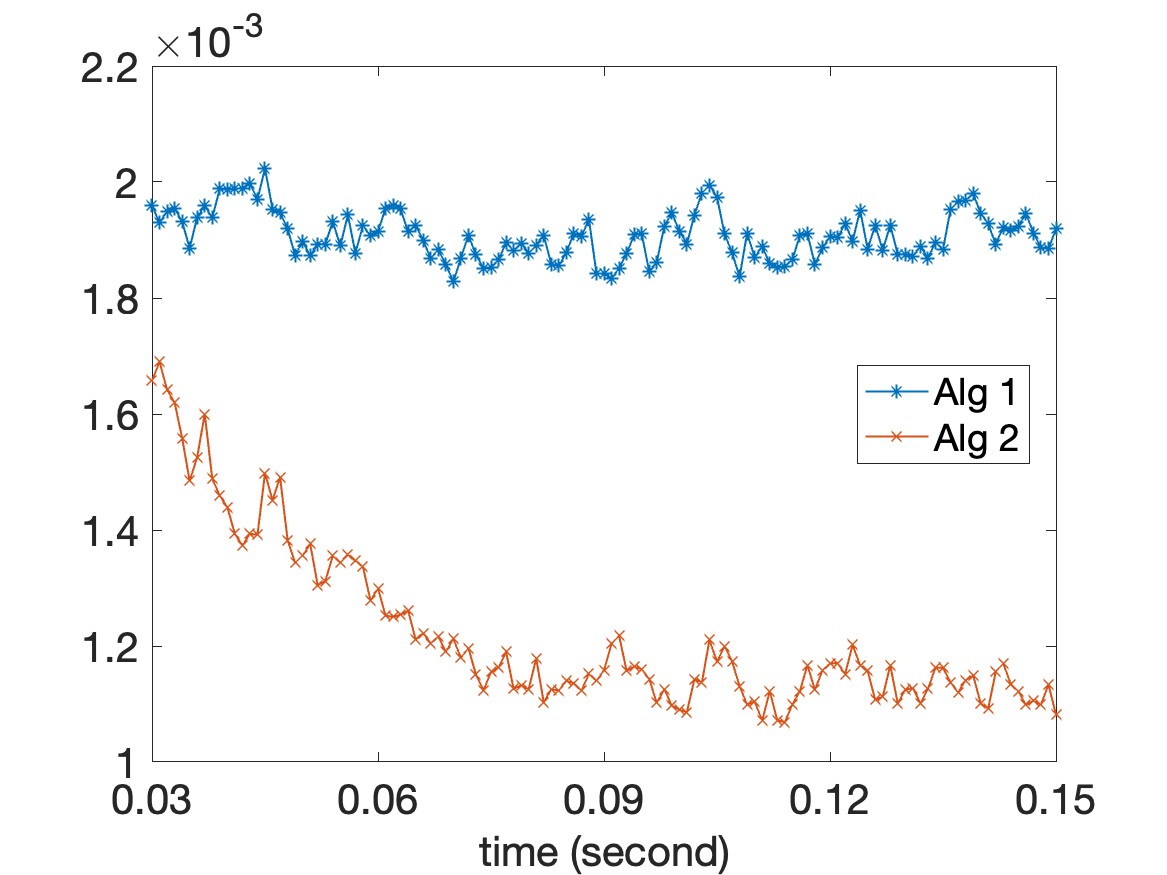}
\caption{Comparison of the relative error in \eqref{eq:err} produced by Algorithm \ref{alg:ETKF} and Algorithm \ref{alg:mETKF} in the time domain $[0.03,0.15]$.}
\label{fig:errorfull}
\end{figure}

Figure \ref{fig:errorfull} displays the relative error given by \eqref{eq:err}
obtained using both Algorithm \ref{alg:ETKF} and Algorithm \ref{alg:mETKF} in the time domain $[0.03,0.15]$. 
The first time interval is omitted since  the solutions have not yet reached a stationary state. Once again we observe that Algorithm \ref{alg:mETKF} yields a smaller posterior error than Algorithm \ref{alg:ETKF} does. 

\subsection{Sparsely sampled observations}
\label{subsection:numerical sparse}
We now present results for the sparsely sampled observation case, that is, where data are observed at every other grid point in the spatial domain. In contrast to the densely sampled case, where we used $\cT = I$ in the cosntruction of $W$, here the localization matrix $\cT$ is tridiagonal with nonzero entries of $1$. We compare results obtained using (1) the localized covariance as described in Algorithm \ref{alg:ETKF}, (2) the localized gradient second moment without clustering (Algorithm \ref{alg:mETKF} with $W$ defined by \eqref{eq: new weighting matrix}), and (3) the localized gradient second moment with clustering (Algorithm \ref{alg:mETKF} with $W$ defined by \eqref{eq: new weighting matrix 2}). 
We correspondingly choose parameters $\alpha = 1.3$ in Algorithm \ref{alg:ETKF} and $\beta$ so that $\max_{i,i'} (W_{i,i'}) = 0.0027$, $1 \le i,i' \le n$ in Algorithm \ref{alg:mETKF}. We furthermore use $dist=1$ in \eqref{eq: dist} when employing clustering in Algorithm \ref{alg:mETKF}.

\begin{figure}[h!]
\centering
\includegraphics[width=0.32\textwidth]{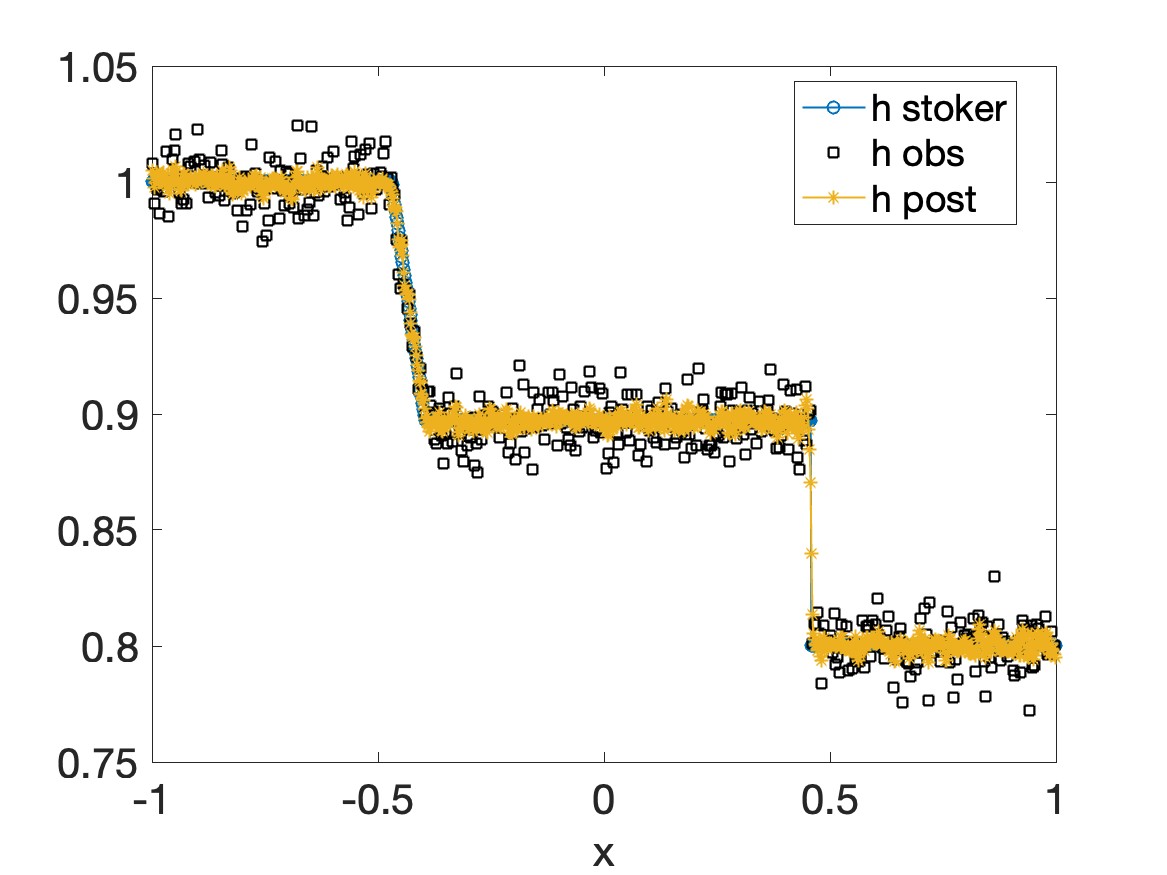}
\includegraphics[width=0.32\textwidth]{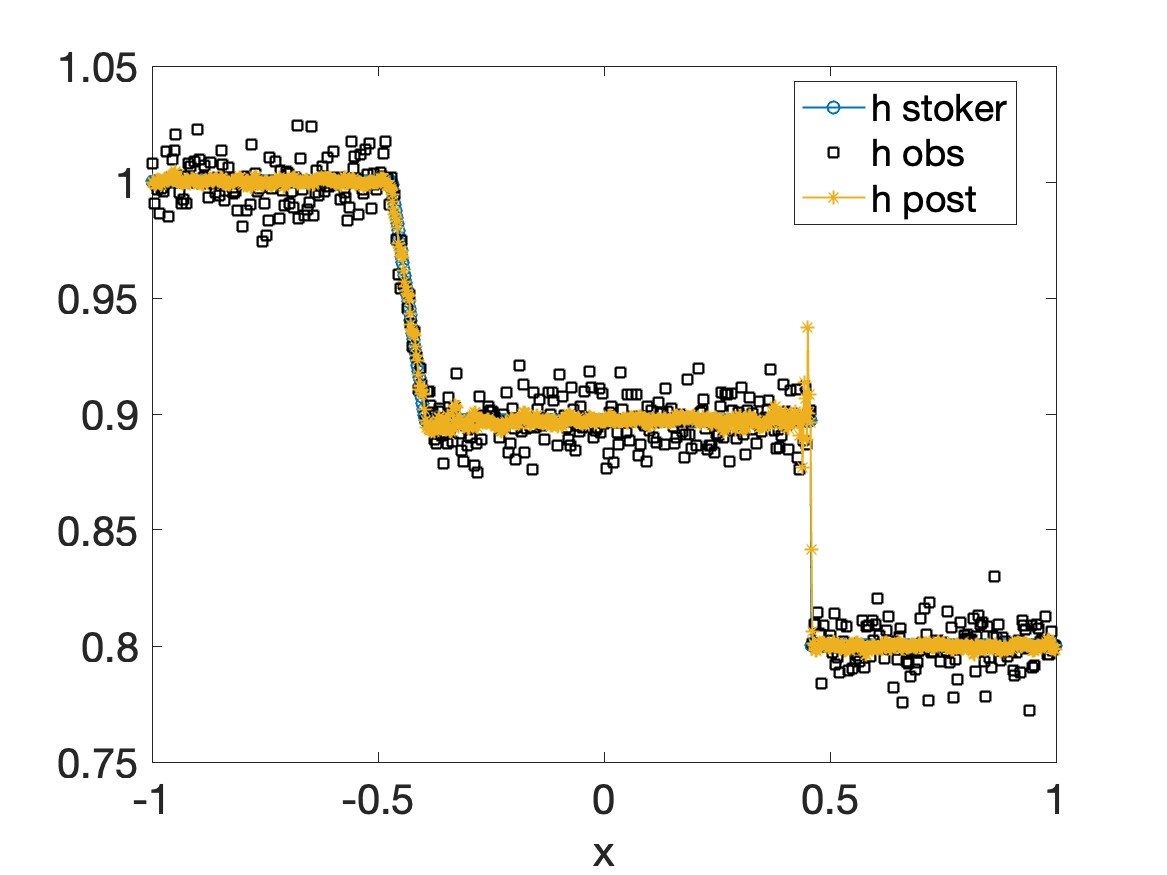}
\includegraphics[width=0.32\textwidth]{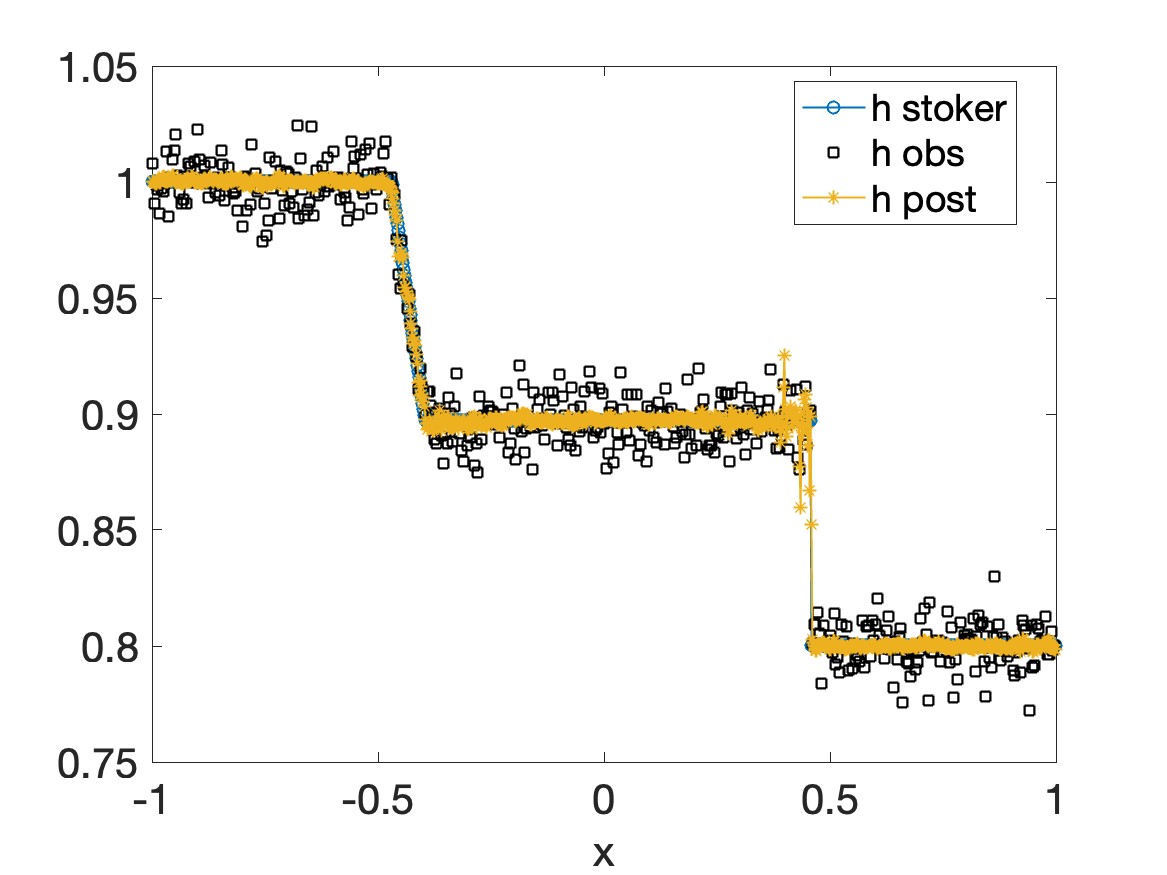}\\
\includegraphics[width=0.32\textwidth]{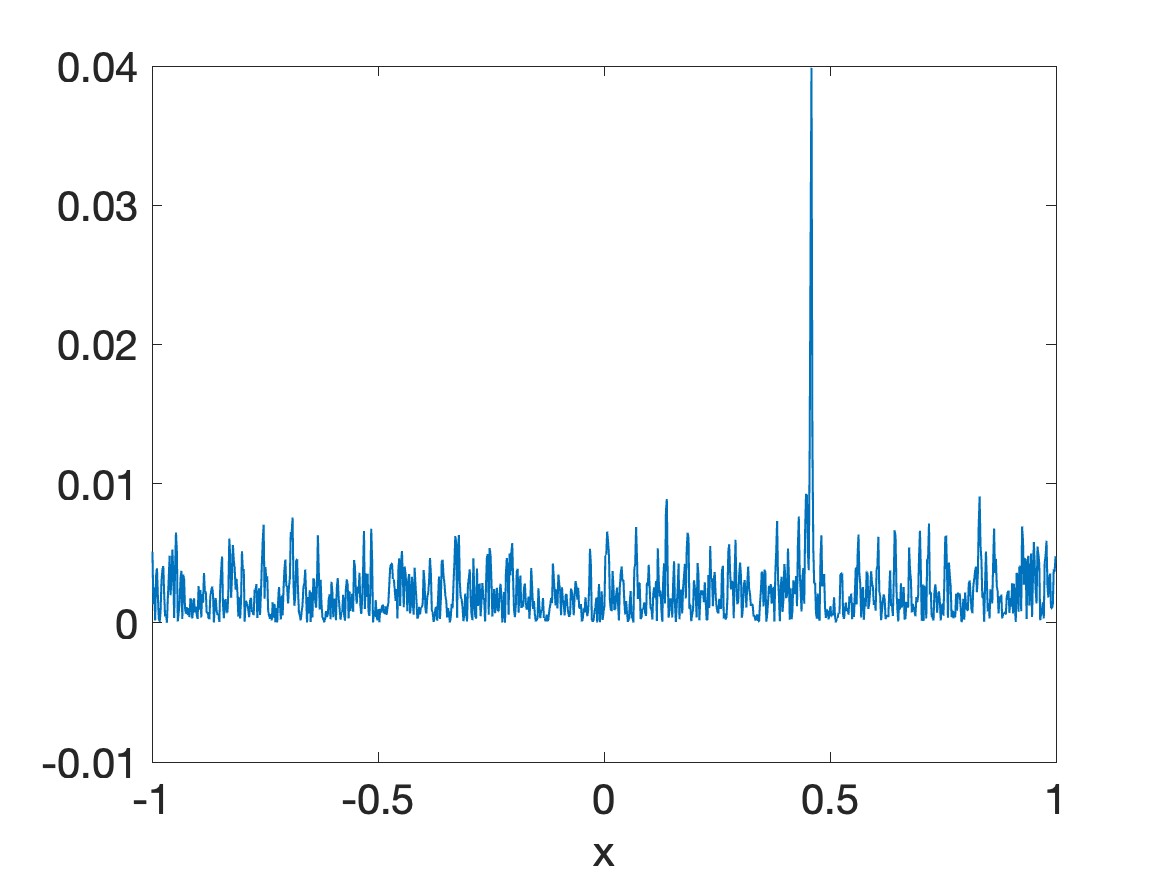}
\includegraphics[width=0.32\textwidth]{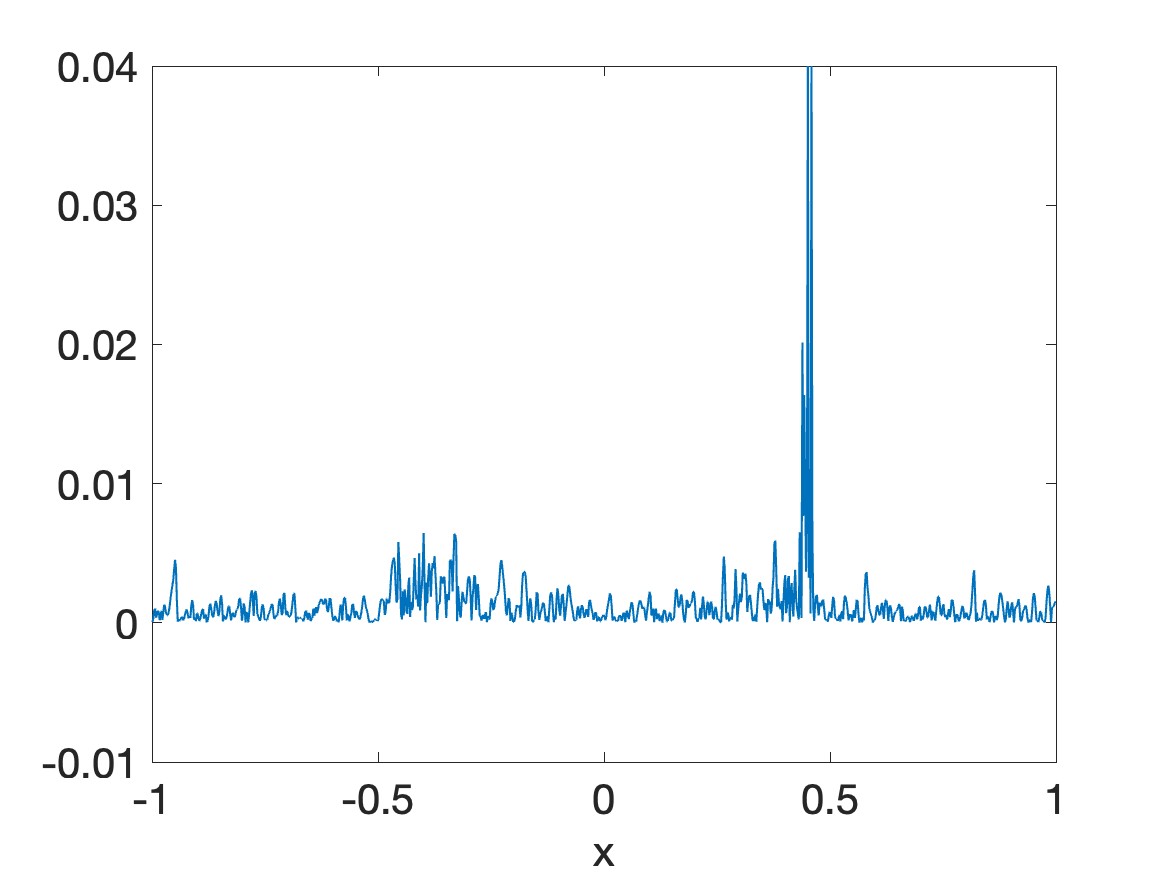}
\includegraphics[width=0.32\textwidth]{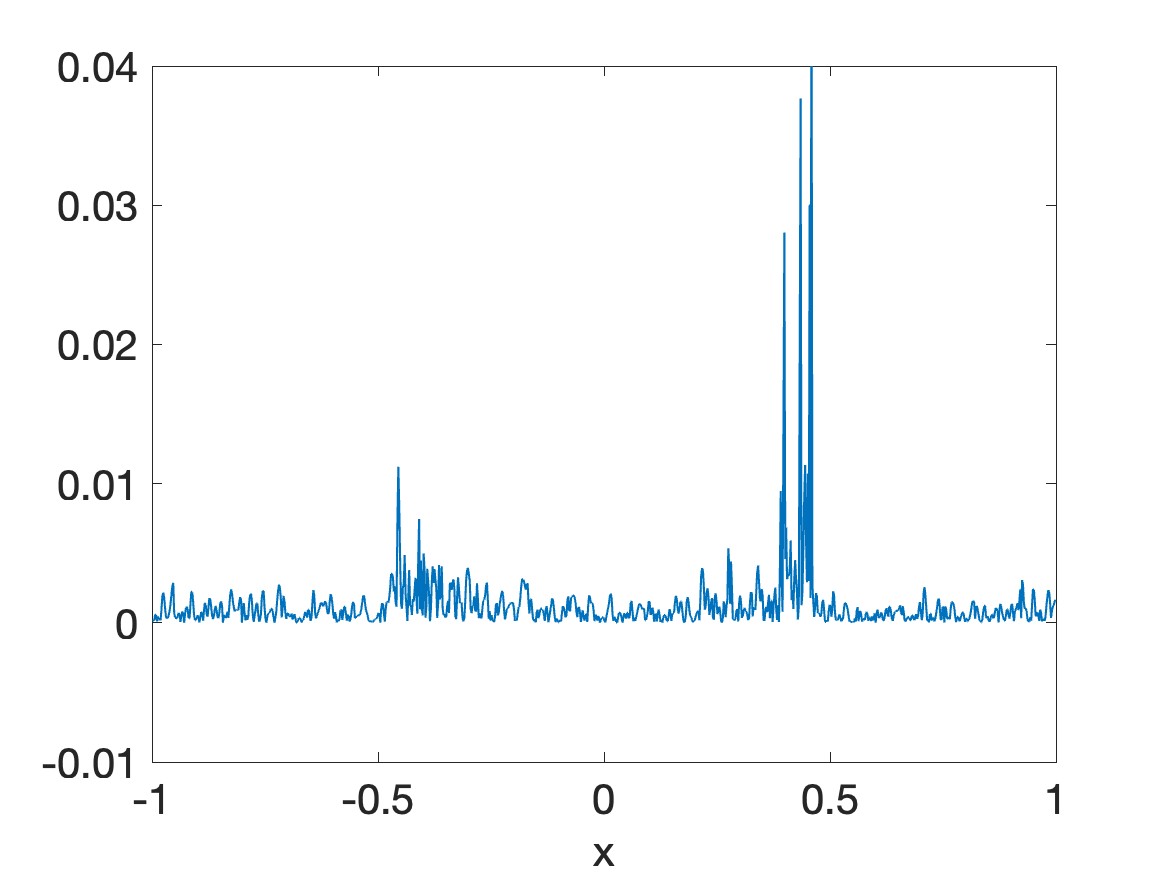}
\caption{(top) Numerical solutions for the posterior mean of depth $h(x,0.15)$ as obtained by (left) Algorithm \ref{alg:ETKF}, (middle) Algorithm \ref{alg:mETKF} without clustering and (right) Algorithm \ref{alg:mETKF} with clustering. (bottom)  Corresponding pointwise posterior error defined by \eqref{eq:error}. }
\label{Fig: sparse sol_}
\end{figure}

Figure \ref{Fig: sparse sol_}(top) shows the numerical solutions of depth $h(x,0.15)$, along with the Stoker solution and the observational data, as obtained by  (left) Algorithm \ref{alg:ETKF}, (middle) Algorithm \ref{alg:mETKF} with $W$ in \eqref{eq: new weighting matrix} and (right) Algorithm \ref{alg:mETKF} with $W$ in \eqref{eq: new weighting matrix 2}, each using the parameter values as already prescribed. Figure \ref{Fig: sparse sol_}(bottom) displays the corresponding pointwise posterior errors using \eqref{eq:error}. It is evident that use of the structurally informed weighted matrix \eqref{eq: new weighting matrix} improves the results, while the additional clustering in \eqref{eq: new weighting matrix 2} helps to mitigate the overshoots at the shock location. Moreover, we note that there are no apparent boundary layer effects resulting from our simulation.

\begin{figure}[h!]
\centering
\includegraphics[width=0.32\textwidth]{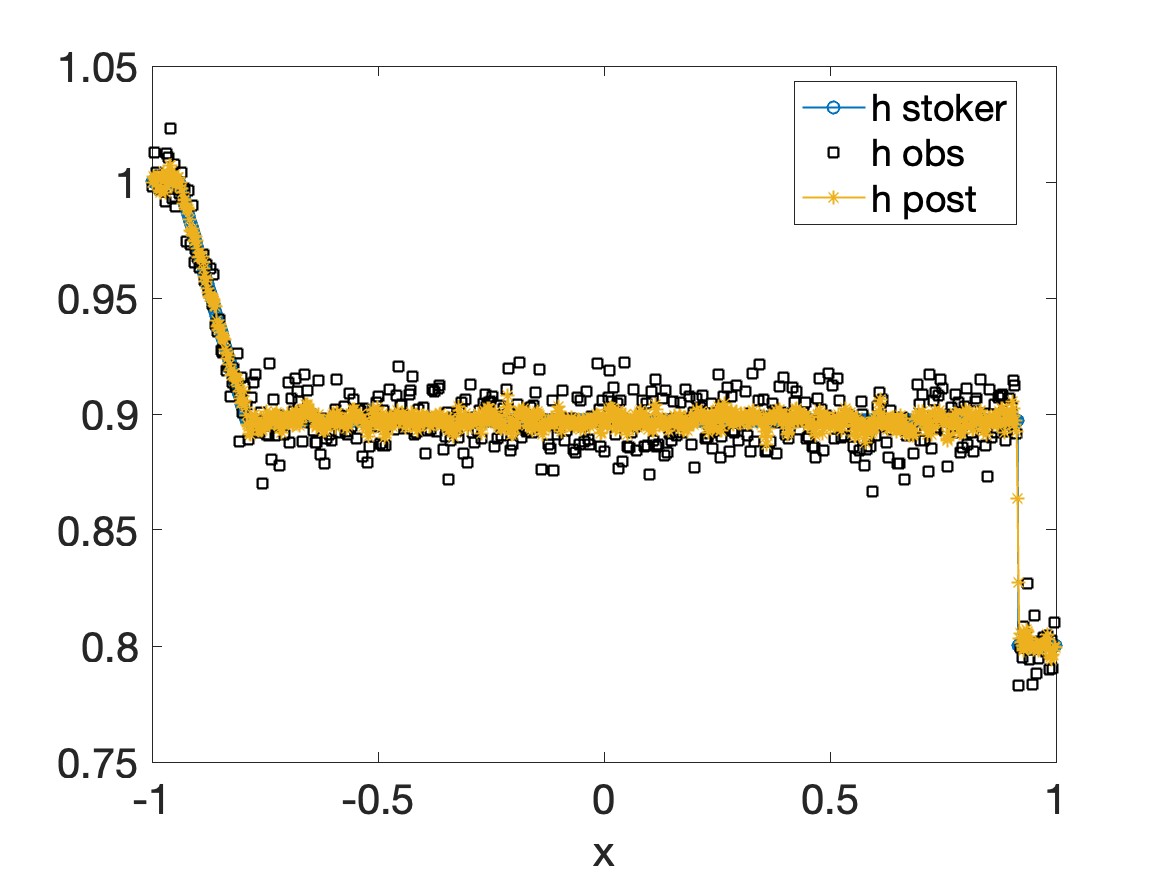}
\includegraphics[width=0.32\textwidth]{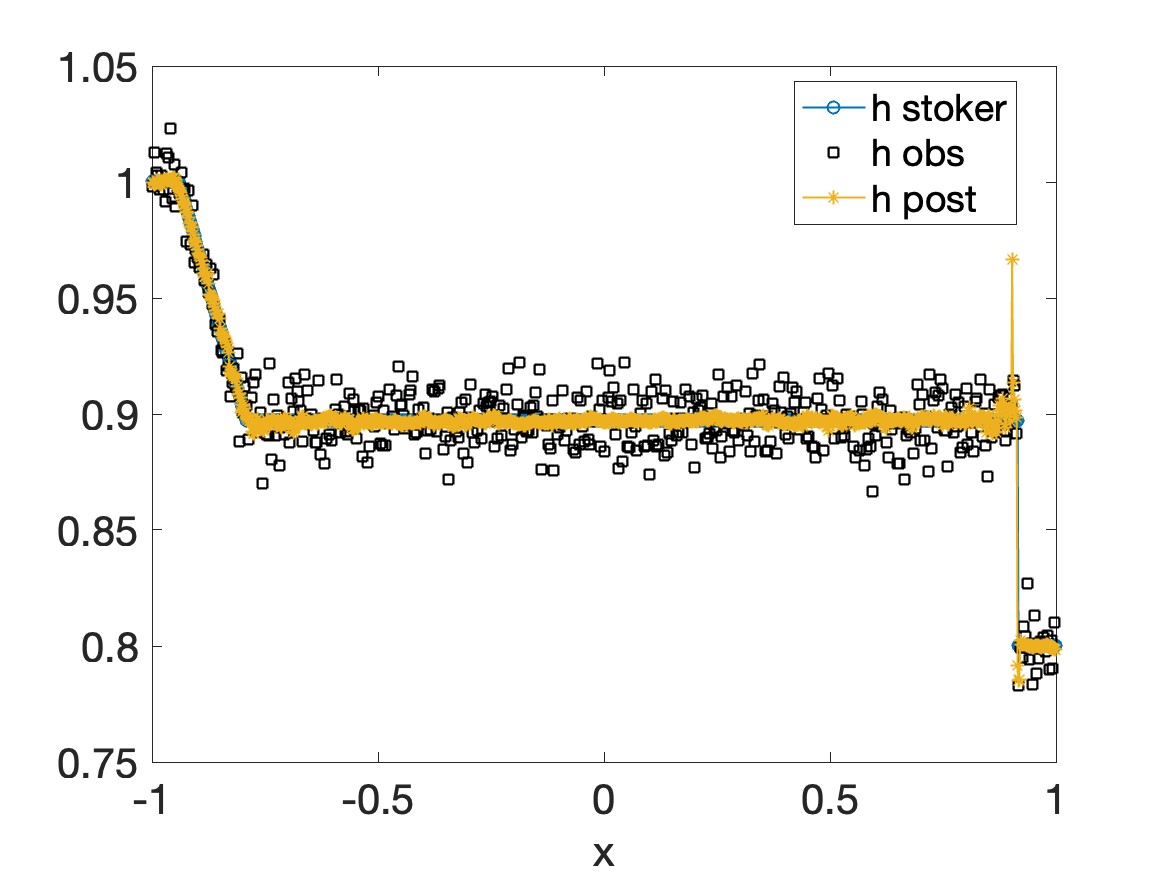}
\includegraphics[width=0.32\textwidth]{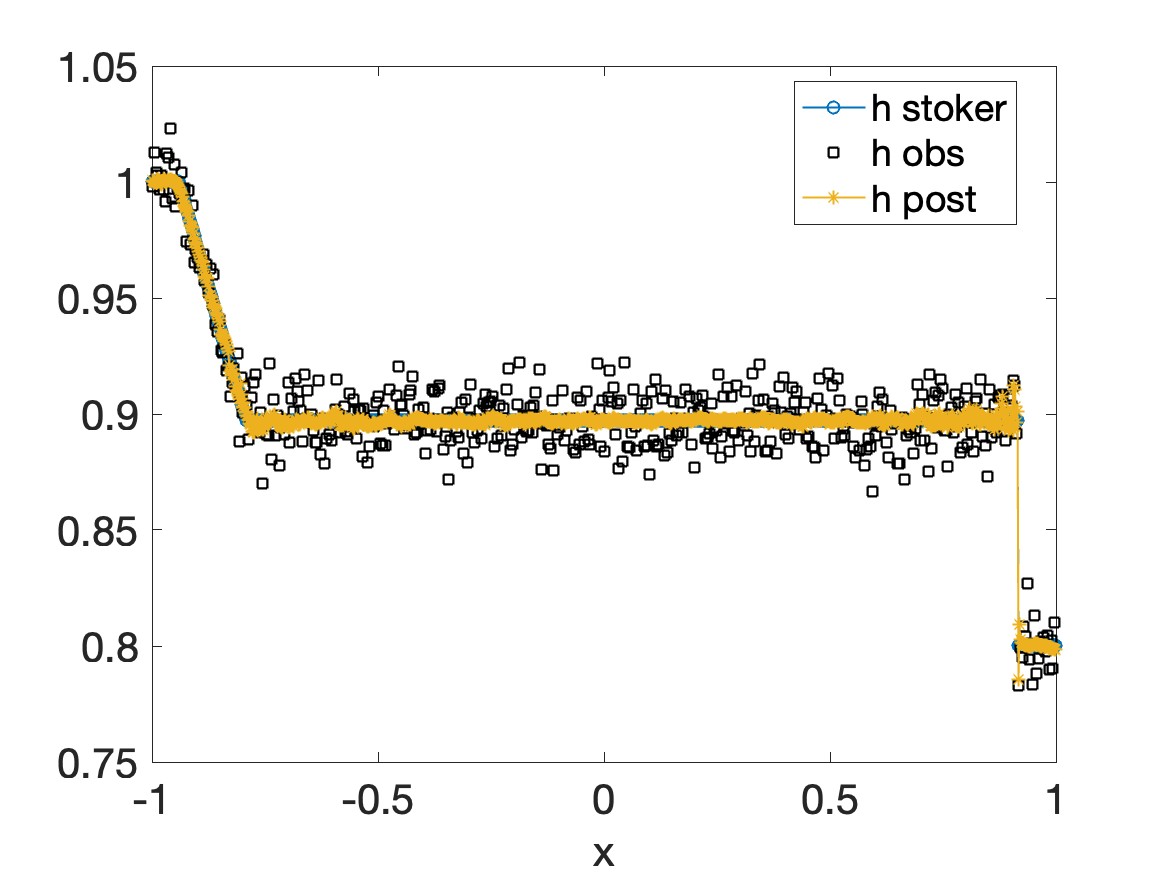}\\
\includegraphics[width=0.32\textwidth]{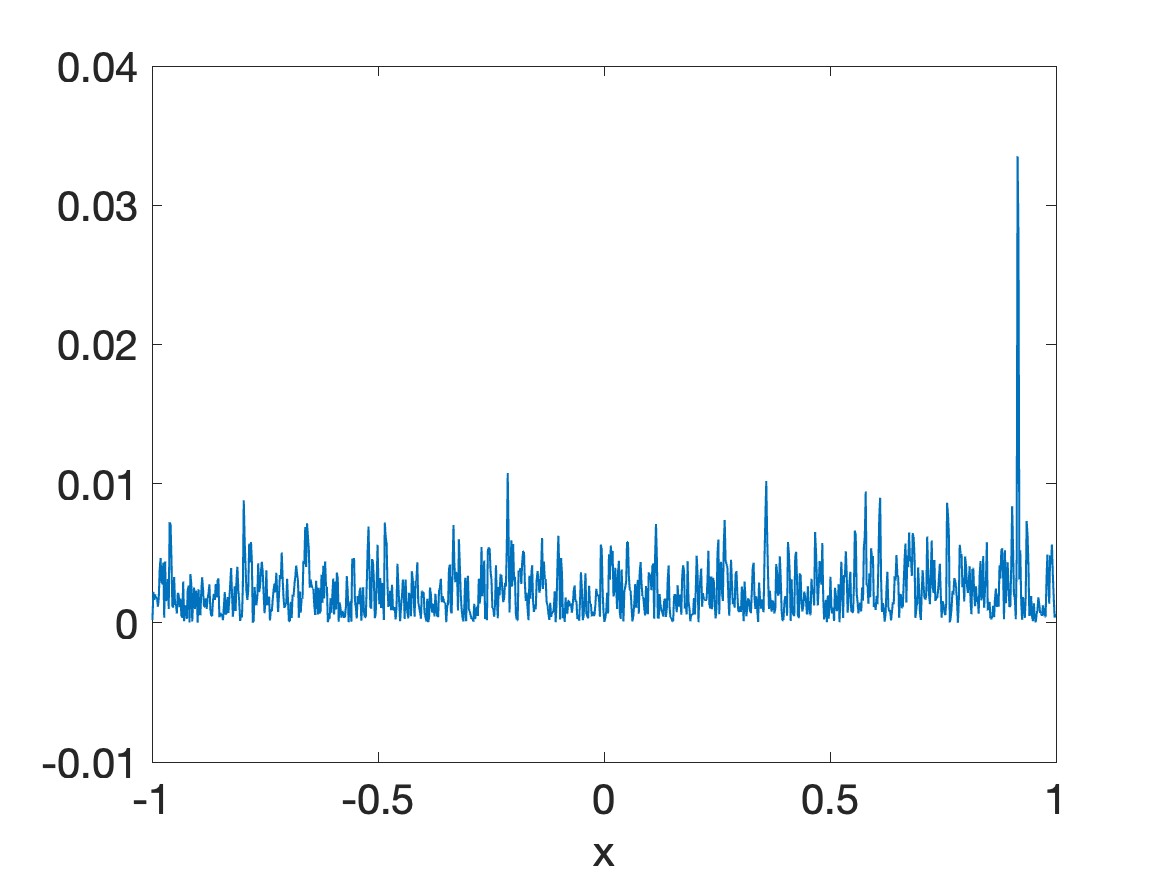}
\includegraphics[width=0.32\textwidth]{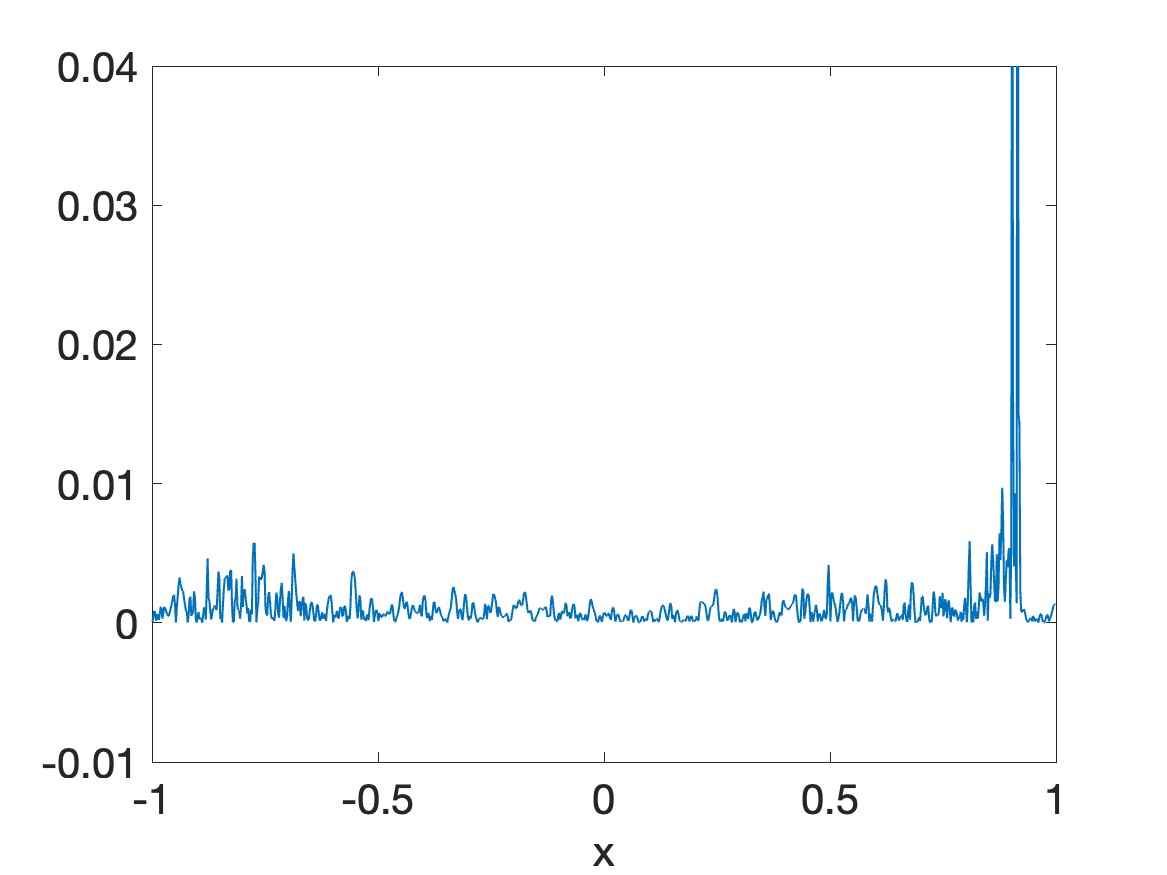}
\includegraphics[width=0.32\textwidth]{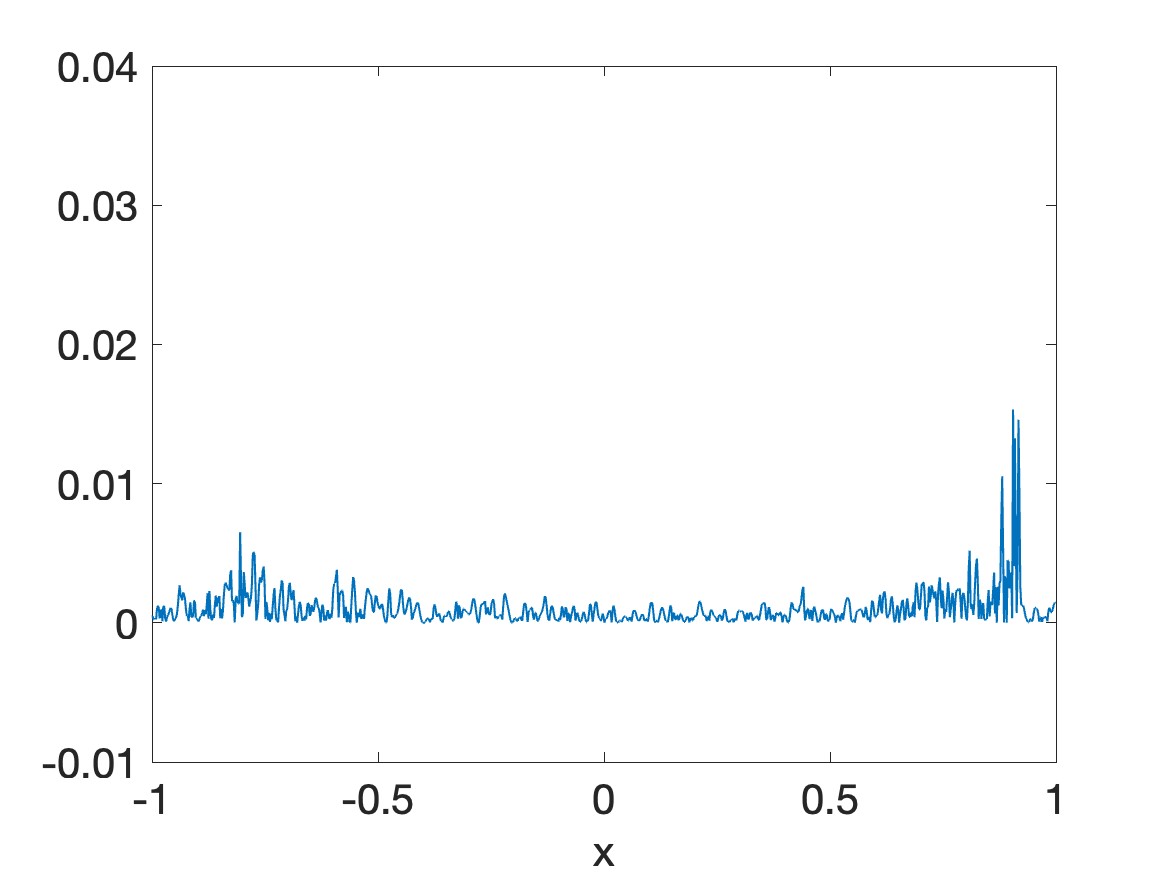}
\caption{(top) Numerical solutions for the posterior mean of depth $h(x,0.3)$ as obtained by (left) Algorithm \ref{alg:ETKF}, (middle) Algorithm \ref{alg:mETKF} without clustering and (right) Algorithm \ref{alg:mETKF} with clustering. (bottom)  Corresponding pointwise squared posterior error defined by \eqref{eq:error}. }
\label{Fig: sparse sol}
\end{figure}
To describe longer term benefits in using our approach for the sparse observation case, Figure \ref{Fig: sparse sol} displays the corresponding numerical solutions and pointwise posterior errors for $h(x,.3)$. Here we observe that while the results appear to be less oscillatory in smooth regions when using $W$ defined by \eqref{eq: new weighting matrix},  incorporating clustering via \eqref{eq: new weighting matrix 2} appears to help mitigate the overshoots at the shock location.  These observations are verified when comparing the magnitudes of the error in Figure \ref{Fig: sparse sol}(bottom).


\begin{figure}[h!]
\centering
\includegraphics[width=0.32\textwidth]{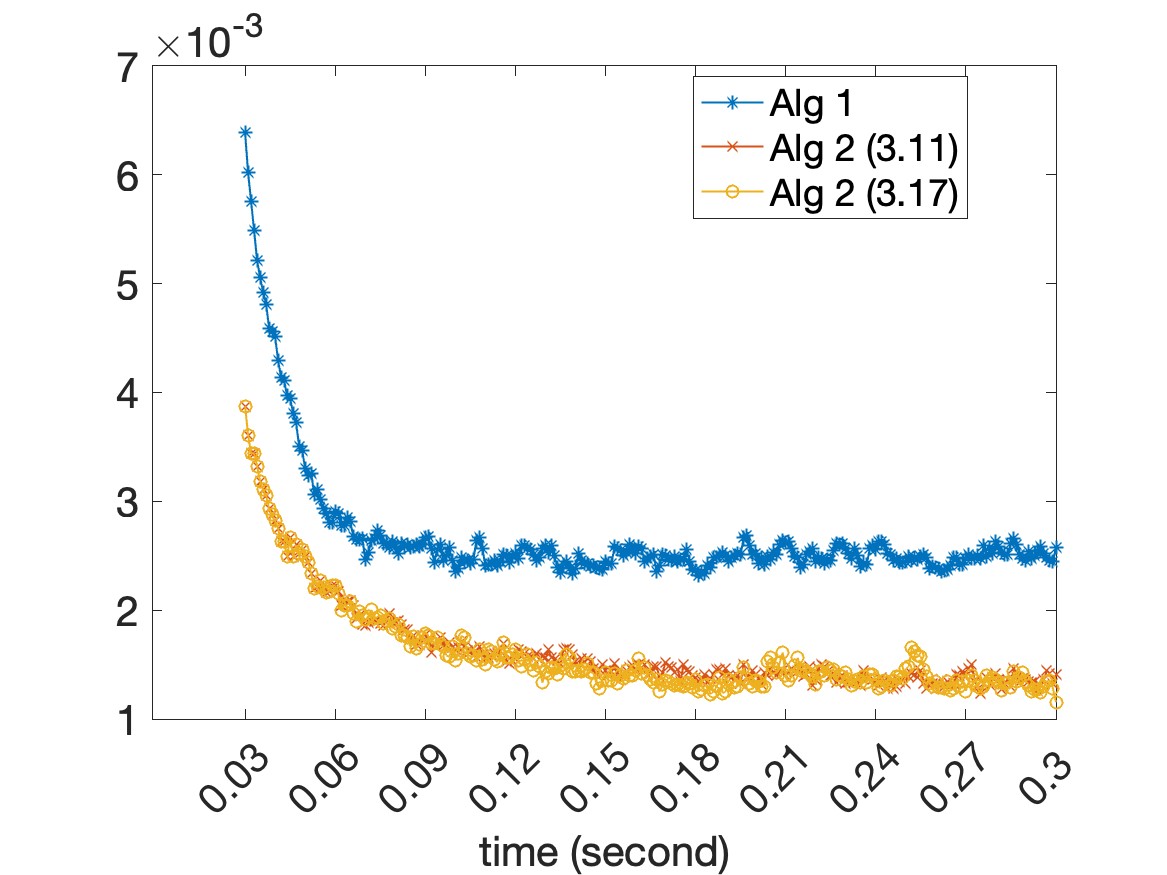}
\includegraphics[width=0.32\textwidth]{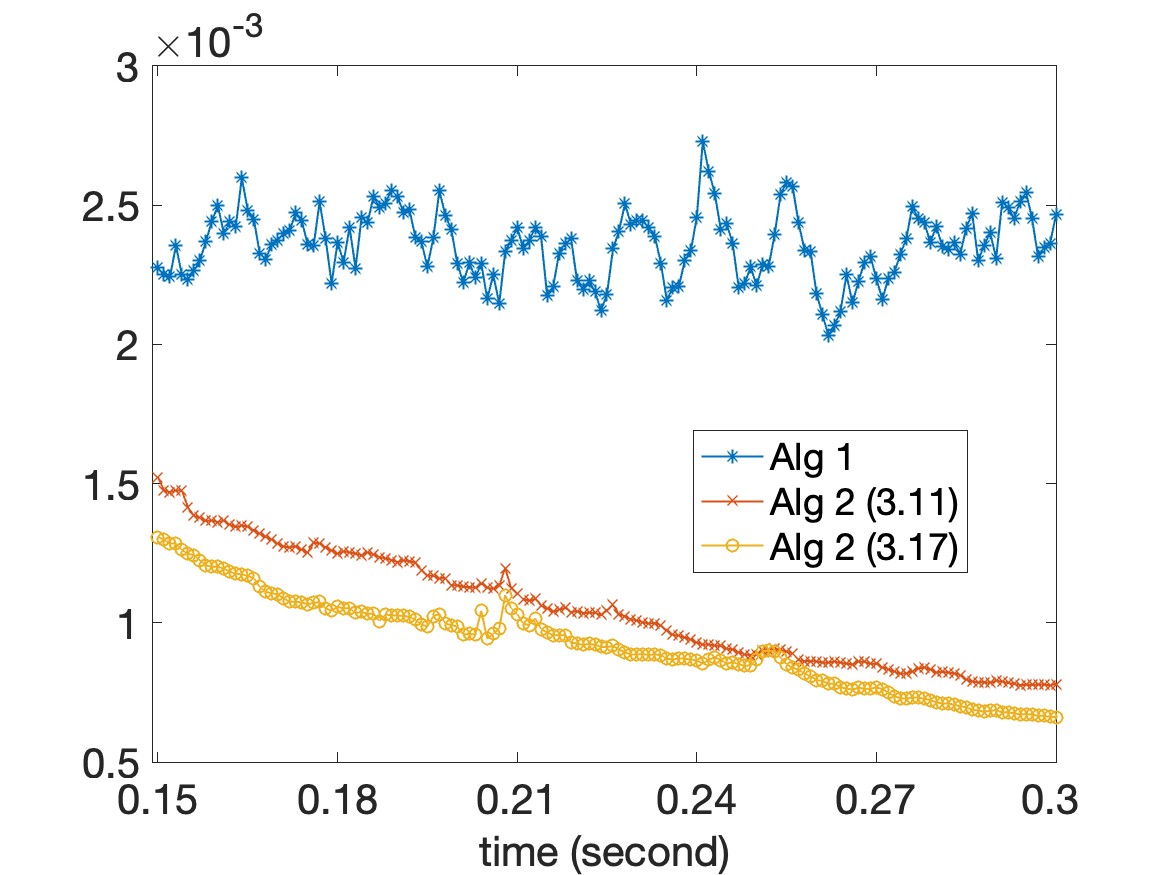}
\caption{{Comparison of Algorithm \ref{alg:ETKF} and Algorithm \ref{alg:mETKF}, both with and without clustering. (Left) relative errors in the time domain $[0.03,0.3]$, (right) relative errors restricted in the spatial domain $[-0.39, 0.39]$ in the time domain $[0.15,0.3]$.}}
\label{fig:errorsparse}
\end{figure}
Figure \ref{fig:errorsparse}(left) displays the relative errors given by \eqref{eq:err} 
obtained using respectively Algorithm \ref{alg:ETKF}, Algorithm \ref{alg:mETKF} with  $W$ in \eqref{eq: new weighting matrix} and Algorithm \ref{alg:mETKF} with $W$ in \eqref{eq: new weighting matrix 2} for the time domain $[0.03,0.3]$. 
The first time interval is omitted since the solutions have not yet reached a stationary state. Figure \ref{fig:errorsparse}(right) shows the relative error modified from \eqref{eq:err} so that the spatial locations are restricted in the smooth domain $[-0.39,0.39]$ for the time domain $[0.15,0.3]$. The results are consistent with those shown in Figure \ref{Fig: sparse sol}.

\subsection{Shallow water equations with oscillatory behavior}
\label{subsection:numerical oscillate} 
As a final test we consider \eqref{eq:SWE} with $u(x,0) = 0$ and oscillating initial conditions for $h$ given by
\begin{gather}
h(x,0)=
\begin{cases}
1+0.03\sin(30x), \quad -1<x<-0.5, \\
1, \quad -0.5<x<0,\\
0.8, \quad 0<x<1.
\end{cases} 
\label{eq: init}
\end{gather}

As these oscillations are always present in the solution, this example provides a greater challenge for data assimilation. In this case there is no analytical solution, so we approximate the true solution\footnote{For notational convenience, we also refer to this approximation as the {\em Stoker} solution.}  using a more highly resolved  fifth order WENO scheme, with 
$$\Delta x = 1\times 10^{-4}, \quad \Delta t = 1\times 10^{-5}.$$

\begin{figure}[h!]
\centering
\includegraphics[width=0.32\textwidth]{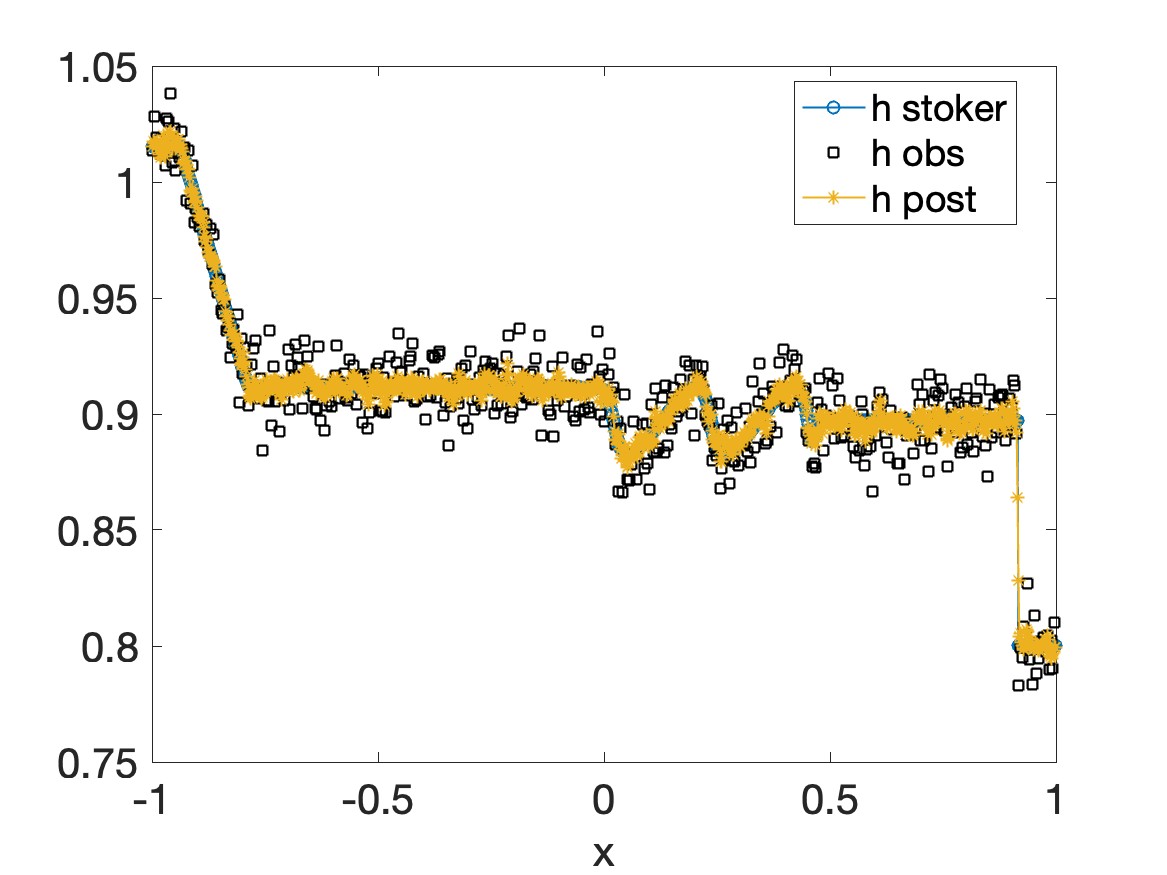}
\includegraphics[width=0.32\textwidth]{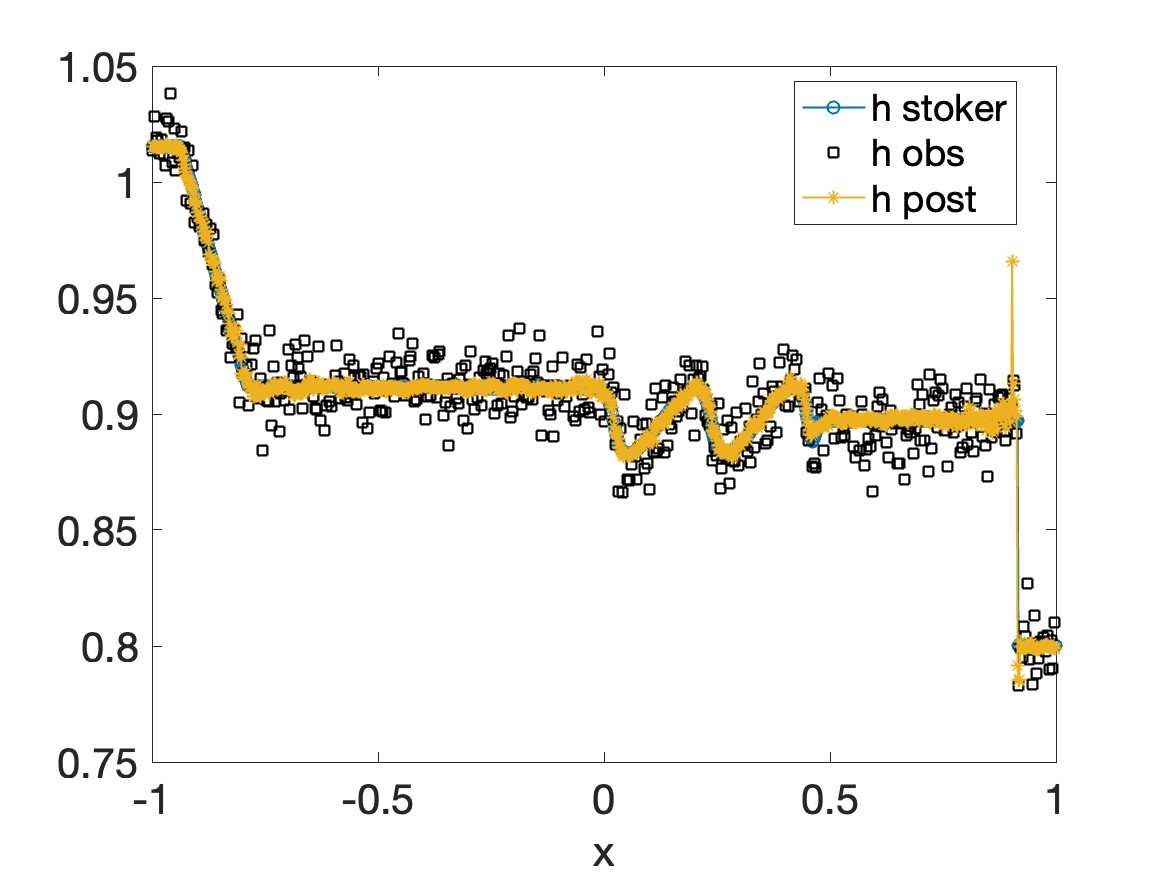}
\includegraphics[width=0.32\textwidth]{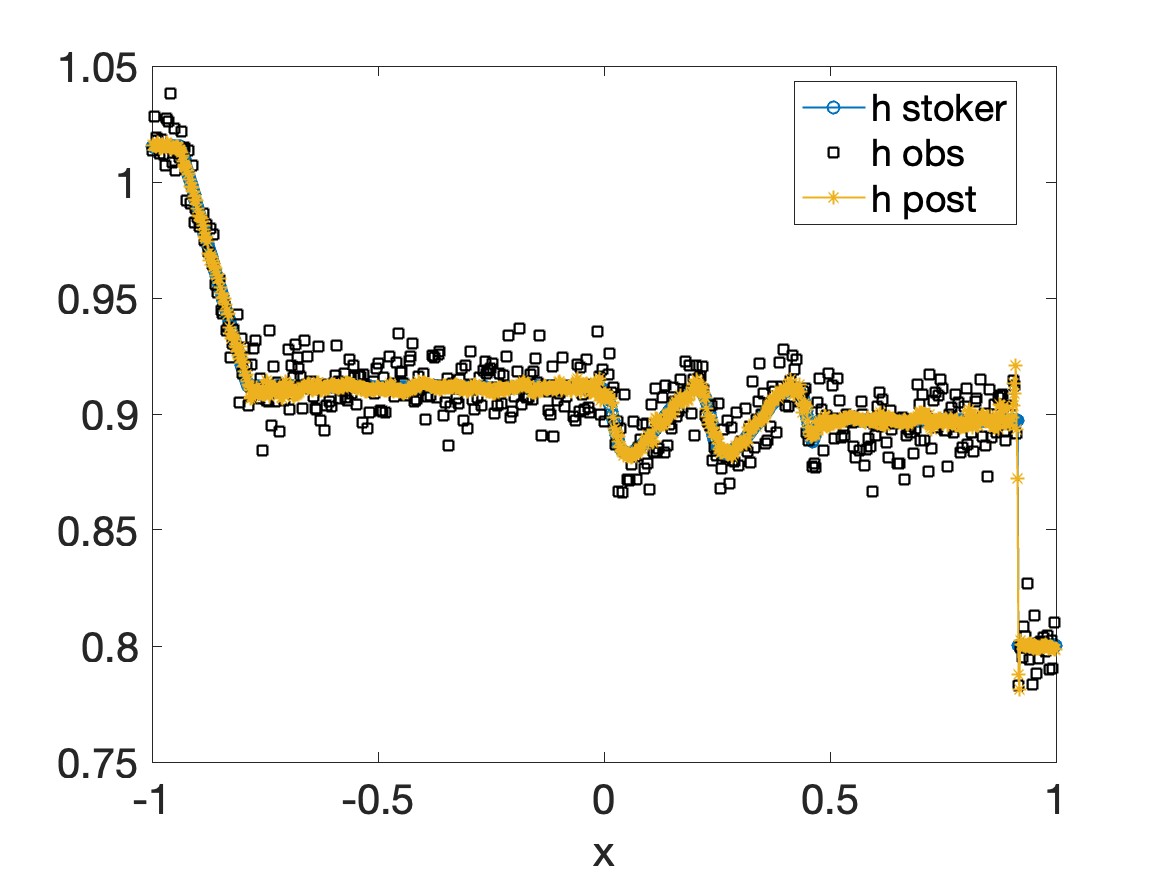}\\
\includegraphics[width=0.32\textwidth]{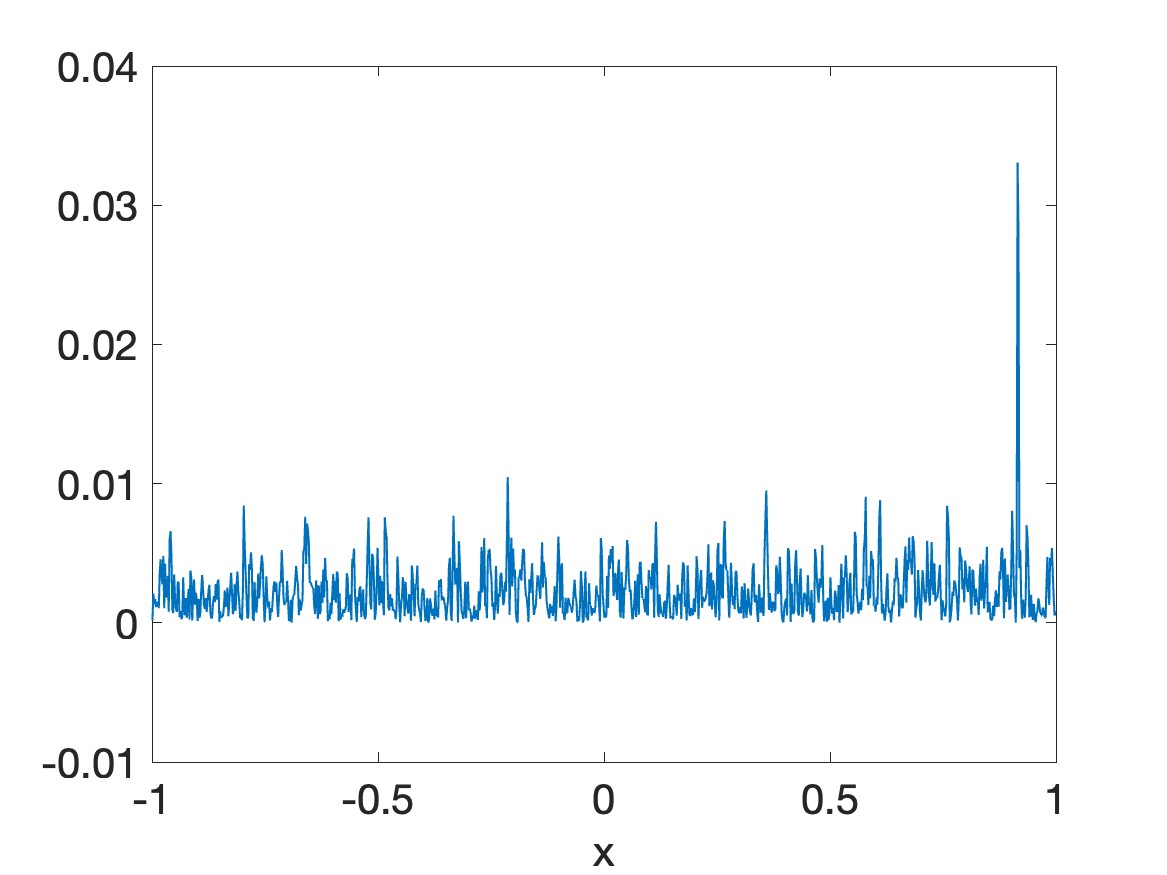}
\includegraphics[width=0.32\textwidth]{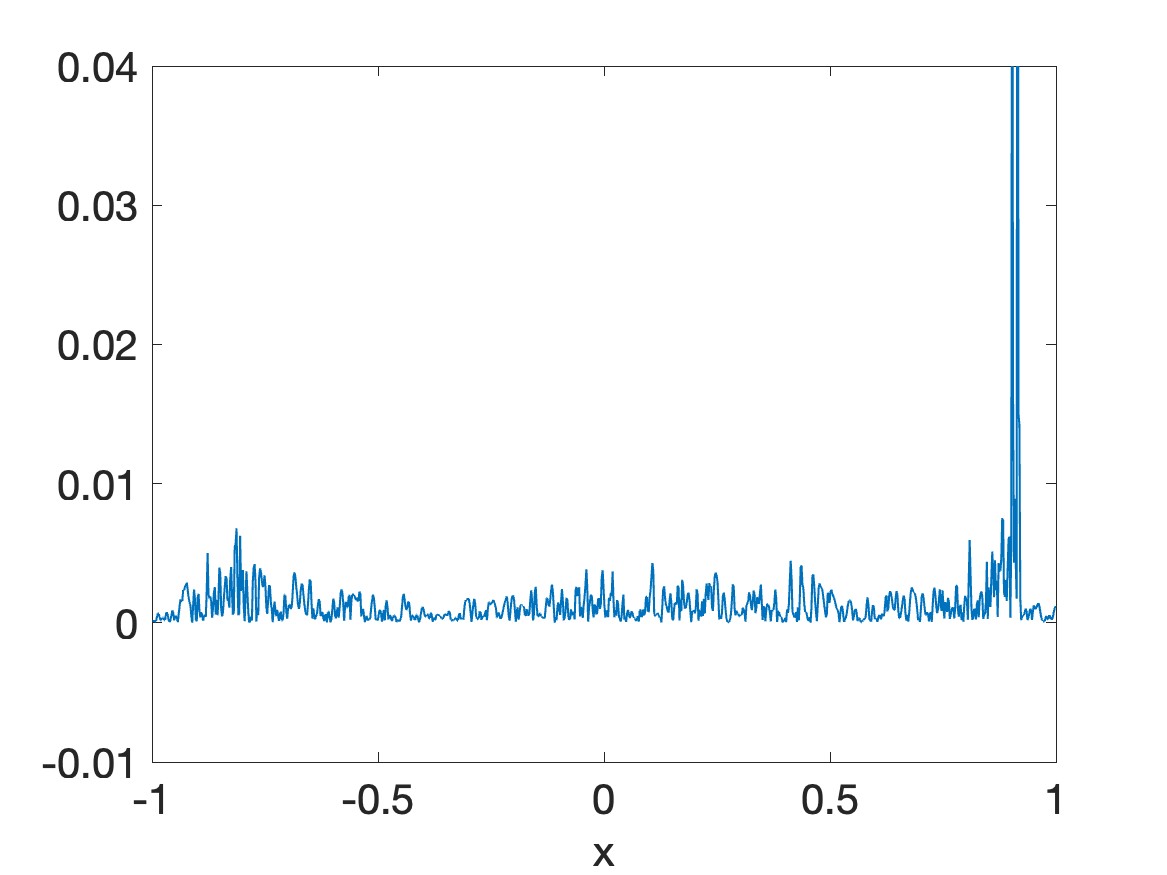}
\includegraphics[width=0.32\textwidth]{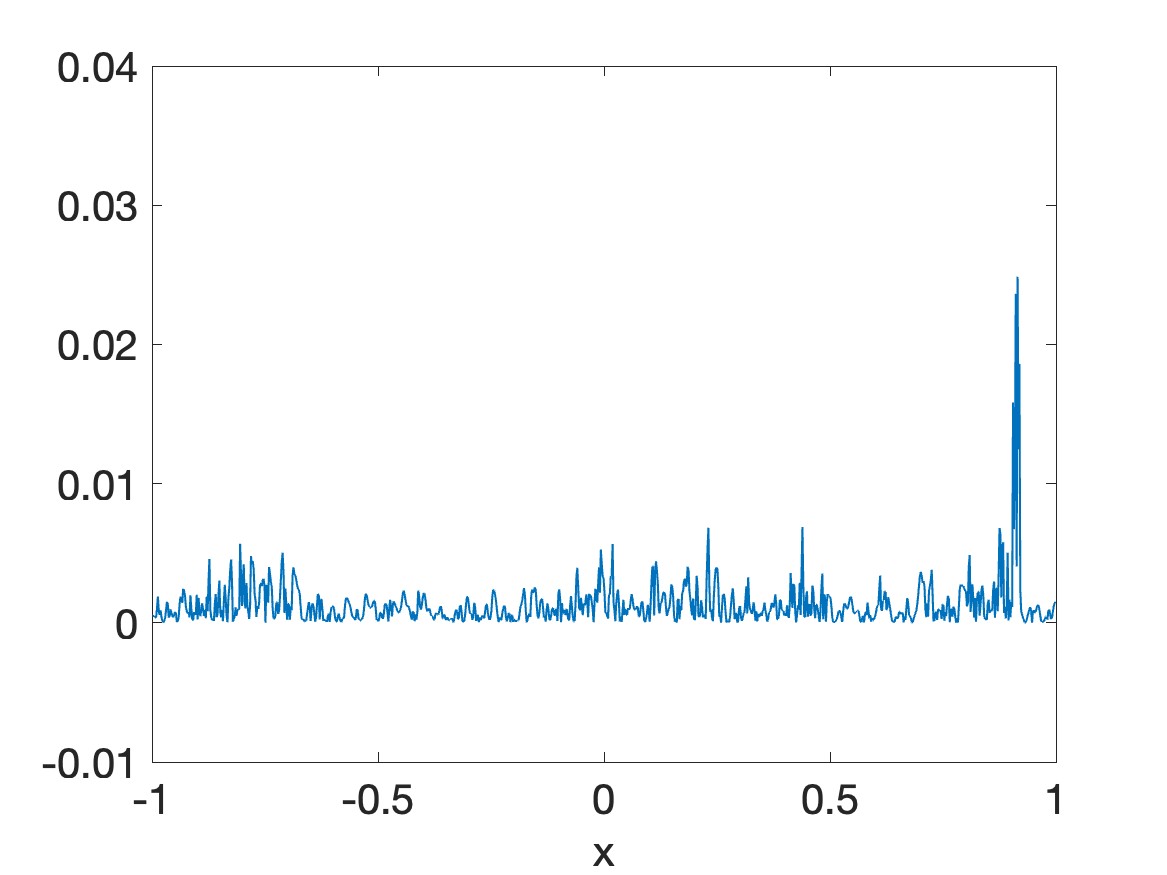}
\caption{(top) Numerical solutions for the posterior mean of depth $h(x,0.3)$ as obtained by (left) Algorithm \ref{alg:ETKF}, (middle) Algorithm \ref{alg:mETKF} without clustering and (right) Algorithm \ref{alg:mETKF} with clustering. (bottom)  Corresponding pointwise posterior error defined by \eqref{eq:error}. }
\label{Fig:oscillate}
\end{figure}

As in Section \ref{subsection:numerical sparse}, we assume we are able to observe data at every other grid point (sparsely sampled  observation case), and choose the same numerical parameters for our experiments. Once again we compare the results using Algorithm \ref{alg:ETKF} with those using Algorithm \ref{alg:mETKF} with and without clustering. 

\begin{figure}[h!]
\centering
\includegraphics[width=0.32\textwidth]{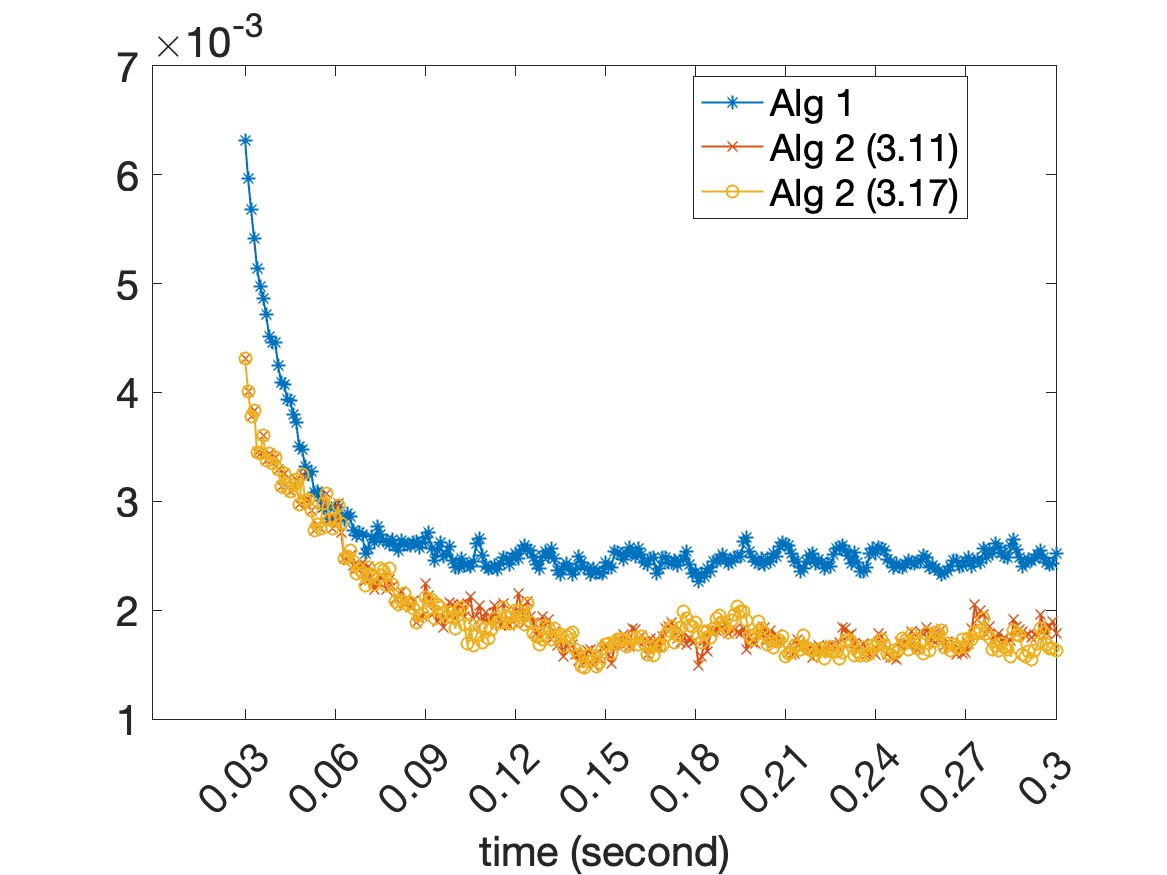}
\includegraphics[width=0.32\textwidth]{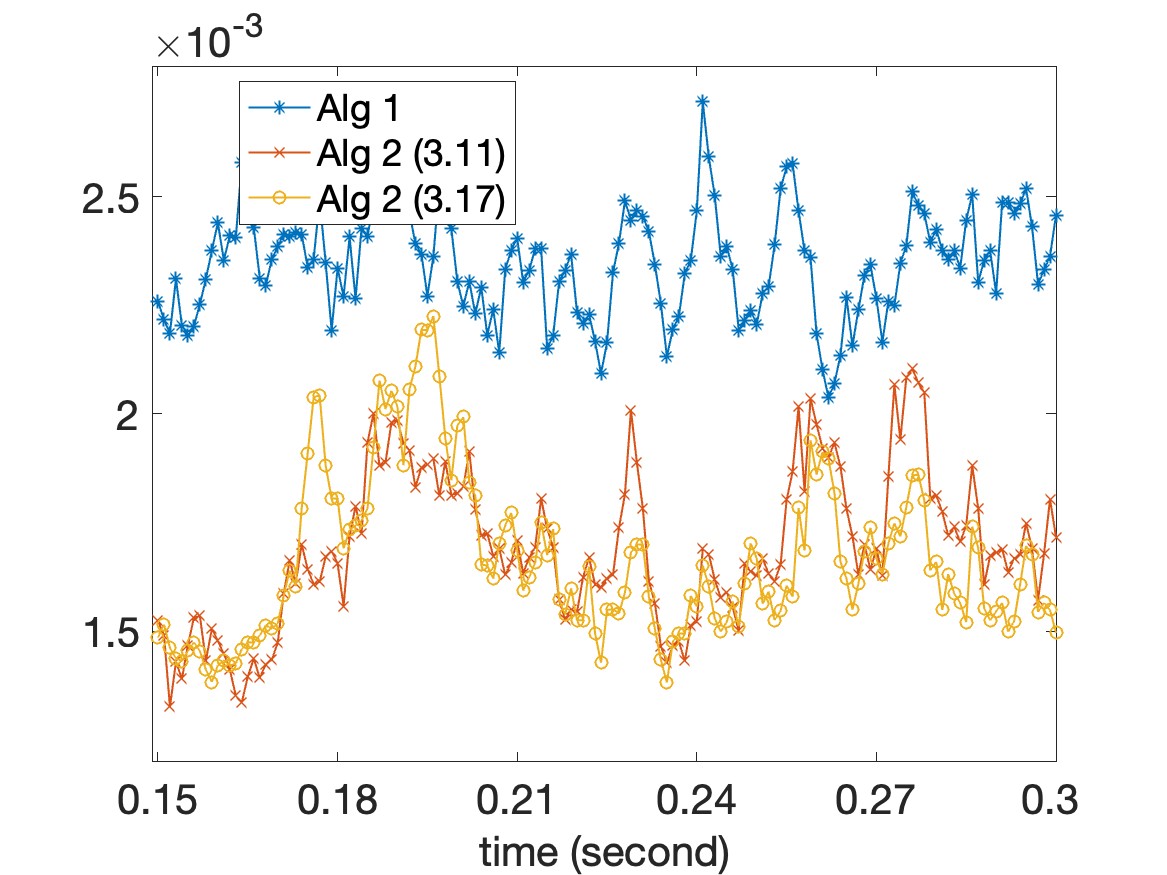}
\caption{
Comparison of Algorithm \ref{alg:ETKF} and Algorithm \ref{alg:mETKF}, both with and without clustering. (Left) relative errors in the time domain $[0.03,0.3]$, (right) relative errors restricted in the spatial domain $[-0.39, 0.39]$ in the time domain $[0.15,0.3]$.}
\label{fig:oscillate_error}
\end{figure}

Figure \ref{Fig:oscillate}(top) compares the posterior mean of the depth $h(x,0.3)$ as obtained by Algorithm \ref{alg:ETKF} and Algorithm \ref{alg:mETKF} with and without clustering, while Figure \ref{Fig:oscillate} (bottom) compares their corresponding pointwise error defined by \eqref{eq:error}. Figure \ref{fig:oscillate_error} shows the relative error for each of these approaches over time, together with the the relative error restricted in the smooth region $[-0.39,0.39]$.  
It is evident that using Algorithm \ref{alg:mETKF} yields more accurate results, and that clustering further helps to reduce the overshoots in the discontinuity regions.

\section{Concluding remarks}
\label{sec:conclusion}


This investigation develops a new structurally informed prior for the statistical data assimilation of non-linear PDEs that contain discontinuities in their solution profiles. Specifically, a new weighting matrix is constructed using second moment gradient information which then replaces the prior covariance matrix typically employed in ensemble-based Kalman filters, resulting in a rebalancing of the data and model terms in the objective function. Our new approach furthermore includes using a numerical method in the prediction step that is able to capture shock information. 
Numerical experiments for both densely and sparsely sampled observations demonstrate that our technique yields more accurate solutions when compared to ETKF, and can be enhanced further when neighborhood information surrounding the local discontinuity structure is also considered (clustering). Our method is robust and requires minimal hand-tuning of parameters. Finally, because the minimization is based on the $\ell_2$ norm, there is a closed-form solution so that no extra optimization algorithms are needed.


Our work presents a new paradigm for incorporating local discontinuity information into the prior design, and there are a number of ideas to explore in terms of both solving specific types of problems as well as practical implementation. For example, here we employed first order finite differencing to compute the gradient.   Higher order methods may be more appropriate if the underlying solution contains more variation or if the spatial grid is less resolved. In some cases this would allow a larger time step and make the prediction step more efficient.  Other shock capturing techniques might also be beneficial in constructing the weighting matrix. For instance, the polynomial annihilation method, \cite{Archibald05,Archibald09,Saxena2009}, is a high order method that is able to detect both shock discontinuities and jumps in the gradient domain in multi-dimensions, and thus may improve the performance of the weighting matrix. Moreover, while $\ell_2$ minimization has the advantage of yielding a closed-form solution, $\ell_1$ regularization methods \cite{CT} are extensively used for sparse signal recovery and generally yield more sharply resolved structures.  Finally, we note that other numerical methods, such as the discontinuous Galerkin (DG) or finite element method, may be more suitable for the prediction step for problems in higher dimensions on non-rectangular grids.  Future work will compare these different approaches for use in our general framework.

\bibliographystyle{abbrv}
\bibliography{mybibfile}

\end{document}